
\documentclass[11pt,a4paper]{article}   
\setlength{\textwidth}{16cm}            
\setlength{\textheight}{24cm}           
\setlength{\oddsidemargin}{-.5mm}       
\setlength{\evensidemargin}{-.5mm}      
\setlength{\topmargin}{-.7cm}	        

\newcommand\leftone{&&\hskip-.5cm}     

\def\LaTeX{%
  \let\Begin\begin
  \let\End\end
  \let\salta\relax
  \let\finqui\relax
  \let\futuro\relax}

\def\UK{\def\our{our}\let\sz s}
\def\USA{\def\our{or}\let\sz z}



\LaTeX

\USA



\salta

\parskip2mm


\usepackage{color}
\usepackage{amsmath}
\usepackage{amsthm}
\usepackage{amssymb}
\usepackage[mathcal]{euscript}




\bibliographystyle{plain}


%

\finqui

\def\Beq{\Begin{equation}}
\def\Eeq{\End{equation}}
\def\Bsist{\Begin{eqnarray}}
\def\Esist{\End{eqnarray}}

\def\Bthm{\Begin{theorem}}
\def\Ethm{\End{theorem}}
\def\Blem{\Begin{lemma}}
\def\Elem{\End{lemma}}
\def\Bprop{\Begin{proposition}}
\def\Eprop{\End{proposition}}
\def\Bcor{\Begin{corollary}}
\def\Ecor{\End{corollary}}
\def\Brem{\Begin{remark}\rm}
\def\Erem{\End{remark}}

\def\Bnot{\Begin{notation}\rm}
\def\Enot{\End{notation}}

\def\Bdim{\Begin{proof}}
\def\Edim{\End{proof}}
\let\non\nonumber




\def\step #1 \par{\medskip\noindent{\bf #1.}\quad}


\def\Lip{Lip\-schitz}
\def\Holder{H\"older}
\def\aand{\quad\hbox{and}\quad}
\def\lsc{lower semicontinuous}

\def\wk{well-known}
\def\lhs{left-hand side}
\def\rhs{right-hand side}


\def\normaliz{normali\sz}

\def\regulariz{regulari\sz}


\def\somma #1#2#3{\sum_{#1=#2}^{#3}}
\def\tonde #1{\left(#1\right)}

\def\graffe #1{\mathopen\{#1\mathclose\}}

\def\<#1>{\mathopen\langle #1\mathclose\rangle}

\def\iot {\int_0^t}
\def\ioT {\int_0^T}
\def\intQt{\int_{Q_t}}
\def\intQ{\int_Q}
\def\iO{\int_\Omega}

\def\dt{\partial_t}
\def\dn{\partial_n}

\def\Vz{V_0}
\def\Vp{V_0^*}

\def\ustar{u^*}

\def\cpto{\,\cdot\,}

\def\checkmmode #1{\relax\ifmmode\hbox{#1}\else{#1}\fi}
\def\aeO{\checkmmode{a.e.\ in~$\Omega$}}
\def\aeQ{\checkmmode{a.e.\ in~$Q$}}

\def\aeS{\checkmmode{a.e.\ on~$\Gamma\times(0,T)$}}

\def\aaQ{\checkmmode{for a.a.~$(x,t)\in Q$}}
\def\aat{\checkmmode{for a.a.~$t\in(0,T)$}}


\def\erre{{\mathbb{R}}}


\def\multibold #1{\def\arg{#1}%
  \ifx\arg\pto \let\next\relax
  \else
  \def\next{\expandafter
    \def\csname #1#1#1\endcsname{{\bf #1}}%
    \multibold}%
  \fi \next}

\def\pto{.}

\def\multical #1{\def\arg{#1}%
  \ifx\arg\pto \let\next\relax
  \else
  \def\next{\expandafter
    \def\csname cal#1\endcsname{{\cal #1}}%
    \multical}%
  \fi \next}


\def\multimathop #1 {\def\arg{#1}%
  \ifx\arg\pto \let\next\relax
  \else
  \def\next{\expandafter
    \def\csname #1\endcsname{\mathop{\rm #1}}%
    \multimathop}%
  \fi \next}

\multibold
qwertyuiopasdfghjklzxcvbnmQWERTYUIOPASDFGHJKLZXCVBNM.

\multical
QWERTYUIOPASDFGHJKLZXCVBNM.

\def\supess{\mathop{\rm sup{\,}ess}}
\def\infess{\mathop{\rm inf{\,}ess}}

\multimathop
dist div dom lip meas sign  supp .


\def\accorpa #1#2{\eqref{#1}--\eqref{#2}}
\def\Accorpa #1#2 #3 {\gdef #1{\eqref{#2}--\eqref{#3}}%
  \wlog{}\wlog{\string #1 -> #2 - #3}\wlog{}}


\def\norma #1{\mathopen \| #1\mathclose \|}
\def\normaV #1{\norma{#1}_V}
\def\normaH #1{\norma{#1}_H}
\def\normaVp #1{\norma{#1}_*}




\def\genspazio #1#2#3#4#5{#1^{#2}(#5,#4;#3)}
\def\spazio #1#2#3{\genspazio {#1}{#2}{#3}T0}

\def\L {\spazio L}
\def\H {\spazio H}
\def\W {\spazio W}

\def\C #1#2{C^{#1}([0,T];#2)}


\def\Lx #1{L^{#1}(\Omega)}
\def\Hx #1{H^{#1}(\Omega)}
\def\Wx #1{W^{#1}(\Omega)}

\def\HxG #1{H^{#1}(\Gamma)}

\def\Luno{\Lx 1}
\def\Ldue{\Lx 2}
\def\Linfty{\Lx\infty}
\def\Huno{\Hx 1}
\def\Hmuno{\Hx {-1}}
\def\Hdue{\Hx 2}
\def\Hunoz{{H^1_0(\Omega)}}


\def\LQ #1{L^{#1}(Q)}


\let\theta\vartheta
\let\eps\varepsilon

\let\TeXchi\chi                         
\newbox\chibox
\setbox0 \hbox{\mathsurround0pt $\TeXchi$}
\setbox\chibox \hbox{\raise\dp0 \box 0 }
\def\chi{\copy\chibox}


\def\soluzione{(\theta,\chi,\xi)}

\def\thetasub_#1{\theta_{\mskip -1.5mu #1}}
\def\soluz #1{(\thetasub_{#1},\chi_{#1},\xi_{#1})}

\def\thetaz{\thetasub_0}

\def\thetaH{\thetasub_{\cal H}}
\def\thetamin{\thetasub_*}
\def\thetamax{\theta^*}
\def\thetaeps{\thetasub_\eps}

\def\thetau{\thetasub_1}
\def\thetad{\thetasub_2}
\def\Istar{[\thetamin,\thetamax]}

\def\lmin{\ell_*}
\def\lmax{\ell^*}
\def\Qep{Q_\eps^+}
\def\Qem{Q_\eps^-}
\def\Oept{\Omega_\eps^+(t)}
\def\Oemt{\Omega_\eps^-(t)}

\def\thetaG{\theta_\Gamma}
\def\thetaH{\theta_{\cal H}}

\def\chiz{\chi_0}

\def\chieps{\chi_\eps}
\def\chimin{\chi_*}
\def\chimax{\chi^*}
\def\chiu{\chi_1}
\def\chid{\chi_2}

\def\xieps{\xi_\eps}

\def\piz{\pi_0}

\let\hat\widehat

\def\Omegatilde{\smash{\widetilde\Omega}}
\def\Qtilde{\smash{\widetilde{\vphantom t\smash Q}}}
\def\Beta{\hat{\vphantom t\smash\beta\mskip2mu}\mskip-1mu}
\def\betatilde{\smash{\widetilde{\vphantom t\smash\beta}}}
\def\betatildeeps{\betatilde_\eps}
\def\betaeps{\beta_\eps}
\def\Betaeps{\Beta_{\!\eps}}
\def\betax{\beta_{\!,x}}
\def\betaepsx{\beta_{\eps,x}}
\def\betat{\beta_{\!,t}}
\def\betaepst{\beta_{\eps,t}}
\def\Betaepst{\Beta_{\!\eps,t}}
\def\betaz{\beta_0}
\def\betau{\beta_1}

\def\cbeta{M_\beta}
\def\deps{\delta_\eps}
\def\zetaeps{\zeta_\eps}
\def\zn{z_n}
\def\betaepsn{\smash{\beta_{\eps_n}}}

\def\lneps{\ln_\eps}
\def\Ln{\mathop{\rm Ln}\nolimits}
\def\Lneps{\Ln_\eps}
\def\rhoeps{\rho_\eps}

\def\Ieps{{\cal L}_\eps}

\def\vn{v_n}

\def\wn{w_n}

\def\Phin{\phi_n}
\def\Psin{\psi_n}
\def\lkeps{\ell_k}
\def\lzeps{\ell_0}
\def\lkpueps{\ell_{k+1}}


\def\cd{c_\delta}
\def\cO{M_\Omega}

\def\cOTm{M_{\Omega,T,m}}
\def\alphastar{\alpha_*}
\def\Cstar{C^*}
\def\lambdap{\lambda_R}

\def\Rie{{\cal R}}





\Begin{document}


{}\vskip-.5cm
\title{Existence and boundedness of solutions\\
for a singular phase field system}

\author{%
\normalsize {\sc Elena Bonetti}$^{(1)}$\\
{\footnotesize e-mail: {\tt elena.bonetti@unipv.it}}\\
\\
\normalsize {\sc Pierluigi Colli}$^{(1)}$\\
{\footnotesize e-mail: {\tt pierluigi.colli@unipv.it}}\\
\\
\normalsize {\sc Mauro Fabrizio}$^{(2)}$\\
{\footnotesize e-mail: {\tt fabrizio@dm.unibo.it}}\\
\\
\normalsize{\sc  Gianni Gilardi}$^{(1)}$\\
{\footnotesize e-mail: {\tt gianni.gilardi@unipv.it}}\\
\\
\\
$^{(1)}${\normalsize Dipartimento di Matematica ``F. Casorati'', Universit\`a di Pavia}\\
{\normalsize Via Ferrata 1, 27100 Pavia, Italy}\\
\\
$^{(2)}${\normalsize Dipartimento di Matematica, Universit\`a di Bologna}\\
{\normalsize Piazza di Porta San Donato 5, 40126 Bologna, Italy}}
\date{}

\maketitle


\vskip-.5cm

\noindent {\bf Abstract.}
This paper is devoted to the mathematical analysis of a 
thermomechanical model describing phase transitions in
terms of the entropy and order structure balance law. 
We consider a macroscopic description 
of the phenomenon and make a presentation of the model.
Then, the initial and boundary value problem is addressed 
for the related PDE system, which contains some nonlinear and 
singular terms with respect to the temperature variable. Existence
of the solution is shown along with the boundedness of both 
phase variable $\chi$ and absolute temperature $\theta$. Finally, 
uniqueness is proved in the framework of a source term 
depending Lipschitz continuously~on~$\theta$.



\noindent{\bf Acknowledgments.}
The authors would like to point out the financial support from the MIUR-COFIN 2006 
research program on ``Free boundary problems, phase transitions 
and models of hysteresis''. The work also benefited from a partial support of
the IMATI of CNR~in~Pavia,~Italy.


\salta

\pagestyle{myheadings}
\markright{}

\finqui


\section{Introduction}
\label{Intro}
\setcounter{equation}{0}

This paper deals with phenomena of phase transitions, of first and second order,
in binary systems (cf., e.g., \cite{F1, FM, G-L}). In a first order phase
transition phenomenon, as in the solid-liquid or liquid-vapor phase change, the
phase transition occurs at a critical temperature, say $\theta_c$: if the
absolute temperature $\theta$ in the body is strictly greater than the
critical temperature $\theta_c$, then the minimum of the energy potential is
attained in one of the pure phases, while if $\theta<\theta_c$ the minimum is
attained in the other phase. In the case when $\theta=\theta_c$ the energy
potential has two minima attained for the two phases, that is phase change may
occur. On the other hand, in the case of second order phase transitions, the
system behaves differently provided $\theta$ is greater or less than the
critical temperature $\theta_c$. Indeed, for high temperatures the energy
potential has only one minimum, while for $\theta<\theta_c$ two minima are
attained with the same values. This second behaviour is characteristic, for
instance, of some solid-solid phase transitions, ferromagnetism, and
superconductivity. 

We are going to investigate a model describing these phenomena by use of
phase-field theories, in terms of temperature $\theta$ and a phase parameter
$\chi,$ that includes the effects of micro-motions and micro-forces responsible
for the phase transition (cf. \cite{fremond} and \cite{gurtin}). Indeed, it is
known that phase transitions are caused by changes occurring at a microscopic
level in the (atomic and/or crystal) structure of the system. These changes are
the effects of micro-forces and motions, which have to be included in the
balance of the energy of the whole system, even if we are providing a
macroscopic description of the phenomenon. We follow the suggestion of
M.~Fr\'emond, who proposed a new balance law on phase transitions 
(cf. \cite{fremond} and \cite{CLSS}). Let us quote a fairly similar 
theory devised by Gurtin~\cite{gurtin} 
for Ginzburg-Landau and Cahn-Hilliard equations.  

Hence, we combine this theory with a model, recently introduced,
based on a {\em reduced energy balance} equation, subsequently termed 
as {\em entropy equation} since entropy is involved. We mainly refer to \cite{BFI}
and especially \cite{bcf} for the derivation of the model and related analytical
results. Note that in \cite{bcf} also thermal memory effects are accounted for
(according to the theory of \cite{gp}), while in the present contribution we are
neglecting them. The main advantage of the model itself is that, once the
problem is solved in a suitable sense, one can obtain directly the positivity of
the temperature, mainly due to the presence of a logarithmic nonlinearity in the
resulting system of partial differential equations. 
This avoids the application of maximum principle
arguments, which are difficult to set in a number of interesting situations.
This fact is pointed out in \cite{bcf}, where the model has been introduced in a
general setting and a global existence result is proved for a weak formulation.
For the sake of completeness, we quote a recent contribution \cite{BFR}, where a
more general (non-smooth) relation between the entropy and the absolute
temperature is considered. Also in this case the physical constraint on the
positivity of the absolute temperature in ensured by the model itself. In
\cite{BFR} the model is recovered by writing the first principle in a dual
formulation (in the sense of the convex analysis) using the entropy in place of
the temperature as state variable. The authors of \cite{BFR} prove existence of
solutions $(\theta,\chi)$ and provide a characterization of the long-time
behaviour of the solution trajectories (see also \cite{BII, BRI}). However,
uniqueness is still an open problem (both in \cite{bcf} and in \cite{BFR}): this
is mainly due to the lack of regularity on the $\theta$ component of the
solution. In \cite{BCFG1, BCFG2} a suitable choice (far from the present
approach) of the heat flux law leads to a linear operator acting on the
temperature, which is of some help in showing existence, uniqueness, and
regularity of the solution. Also the large time behaviour and the $\omega$-limit
set are investigated in \cite{BCFG2}.

\step The model

In the actual contribution, we mainly proceed according to \cite{bcf}. Indeed,
we consider a two-phases system, located in a smooth and bounded domain
$\Omega\subset\erre^3$, and look at its evolution during a finite time interval
$(0,T)$. We denote by $\Gamma$ the boundary $\partial\Omega$.
The thermomechanical equilibrium of the system is described in terms of state
variables and governed by the free energy, while the dynamics ensues from the
presence of a pseudo-potential of dissipation (depending on dissipative
variables). Let us point out that the properties of the pseudo-potential of
dissipation ensure thermodynamical consistency of the model. 

We do not consider mechanical effects, so that the variables of the system are
just the absolute temperature $\theta$ and a phase parameter $\chi$, related to
the proportion of one phase with respect to the other. In general, $\chi$
attains its physical admissible values in a range $[\chi_*,\chi^*]$ (e.g.,
$\chi\in [0,1]$) and this physical constraint has to be ensured by the model
itself. Hence, we derive the equations of the system by thermomechanical laws.
More precisely, we use the approach followed 
by Fr\'emond (yielding the evolution of the phase
parameter), and the first and second principle of Thermodynamics (from which we
recover a balance equation ruling the evolution of the temperature).

We first discuss, with some details, the derivation of the evolution equation
for the phase parameter $\chi$. It is known that the phase transition can be
described as a change in the order structure of the thermomechanical system we
are considering, so that the phase parameter $\chi$ can be interpreted also as
an ``order parameter''. More precisely, below a critical temperature it is
observed that the structure order of many materials is greater than above. For
instance, in the solid-liquid phase transition the solid phase has a greater
order structure due to the crystal symmetry group. Analogously, we may think to
ferromagnetism where the magnetic moments are aligned below the Curie
temperature. Obviously, order change in the structure of the system occurs at a
microscopic level. However, the parameter $\chi$ is a macroscopic parameter
whose evolution is governed by a balance law on the order structure, responsible
for the phase transition at a microscopic level. Thus, the evolution of the
order parameter $\chi$ can be derived from the balance conditions
\begin{equation}\label{eqBH}
B-\div{\bf H}=0 \quad\hbox{in } \, \Omega  \times (0,T) , \ \quad{\bf H}\cdot{\bf n}=0 
\quad\hbox{on } \, \Gamma  \times (0,T) 
\end{equation}
${\bf n}$ being the outward normal to the boundary $\Gamma$, where the scalar 
quantity $B$ can be interpreted as an internal order structure  density per 
units of concentration $\chi$, and ${\bf
H}$  represents an order structure flux vector. 
Note that here we are assuming that no
external action is provided to the system. 

As far as the description of the evolution of the temperature is concerned,
we use a balance law in which higher order dissipative contributions
are neglected by means of the small perturbations assumption
(cf.~\cite{germain}). The resulting equation rules
the evolution of the entropy $s$ of the system in
terms of the entropy flux ${\bf Q}$ and an external source $R$,
possibly depending on the state field, that is
\begin{equation}\label{entropy}
s_t+\hbox{div }{\bf Q}= R \quad\hbox{in } \, \Omega  \times (0,T) .
\end{equation} 
Then, \eqref{entropy} is combined with boundary conditions on the entropy flux,
e.g., if no flux is assumed through the boundary, then one states that ${\bf
Q}\cdot{\bf n}=0$ on $\Gamma  \times (0,T)$. Another possibility, actually
followed by us in the analysis, is to prescribe the value of the temperature on
the boundary.

Referring to \cite{bcf} and {\cite{BFR}, we prefer
not to detail here the derivation of the model.  However, 
for the sake of completeness, we just recall that 
\eqref{entropy} can be obtained, in the framework of small perturbations
assumptions, dividing by $\theta>0$ the energy balance
\begin{equation}\label{energybal}
e_t+\hbox{div }{\bf q}=r+B\chi_t+{\bf H}\cdot\nabla\chi_t.
\end{equation}
In \eqref{energybal} $e$ denotes the internal energy of the system, ${\bf q}= \theta {\bf Q}$ represents the heat flux, $r$ stands for an external source, while $B\chi_t+{\bf
H}\cdot\nabla\chi_t$ is internal and comes from the order structure.
We will add some comments below
on the relation between the two equations 
\eqref{entropy} and~\eqref{energybal}.

\step Free energy

We specify the involved physical quantities with the help of two energy
functionals: the free energy $\Psi$, depending on the state variables and
accounting for the thermo-mechanical equilibrium of the system, and the
pseudo-potential of dissipation $\Phi$ (see \cite{moreau}), defined for
dissipative variables and responsible of the evolution of the system. More
precisely, we consider as state variables the absolute temperature $\theta$, the
order parameter $\chi$, and its gradient $\nabla\chi$.  
It is known from Thermodynamics that the free energy is a concave function
with respect to the temperature, while there are not constraints concerning the
dependence on the other variables. In the present contribution, we choose the
functional $\Psi$ of the following form
\begin{equation}\label{free-en}
\Psi(\theta,\chi,\nabla\chi):=-\frac{c_0}2\theta^2+F(\chi)\theta_c+G(\chi)\theta+\frac\nu 2|\nabla\chi|^2
\end{equation}
the constants $c_0,\, \nu$  being positive and  $\theta_c>0$ representing the
critical value of the temperature for the phase transition.
Note that the purely caloric part in the free energy $- (c_0/2)\theta^2$ is,
in fact,  concave with respect to the temperature. Then, the functions $F$ and $G$ 
characterize the behaviour of the phase transition. For instance, in a first
order phase transition, as vapor-liquid or liquid-solid, they can be prescribed
as follows
\begin{equation}\label{spe1}
F(\chi)=\frac{\chi^4}4-\frac{\chi^3}3,\quad G(\chi)=
\frac{\chi^4}4-\frac{2\chi^3}3+\frac{\chi^2}2
\end{equation}
while in a second order phase transition, as for superconductivity 
or ferromagnetism, $F$ and $G$ can be written as
\begin{equation}\label{spe2}
F(\chi)=\frac{\chi^4}4-\frac{\chi^2}2,\quad G(\chi)=\frac{\chi^2}2.
\end{equation}
Let us remark from the beginning that both the cases \eqref{spe1} and
\eqref{spe2} comply with our general assumptions \accorpa{hpFG}{hpFGbis} provided
you take $\chimin  \leq 0 $ and $\chimax \geq 1$. 

\step Admissible values for the phase variable

The physical constraint on $\chi$, that is $\chi_* \leq \chi\leq \chi^*$, is not
a priori guaranteed by the choice of the free energy \eqref{free-en}: in
particular, the functions $F$ and $G$ prescribed in \eqref{spe1} and
\eqref{spe2} are smooth on the whole of $\erre$ (as in \cite{BFG}). A different
possible choice for $\Psi$ introduced in the literature (cf., e.g.,
\cite{fremond, bcf}) and accounting for this constraint  is
$$
\Psi(\theta,\chi,\nabla\chi)=-\frac{c_0}2\theta^2+I_{[\chi_*,\chi^*]}(\chi)+{\cal 
F}(\chi)+\theta{\cal G}(\chi)+\frac\nu 2|\nabla\chi|^2
$$
where ${\cal F}, \, {\cal G}$ are sufficiently smooth functions characterizing
the phase transition.
In this case the free energy is defined for any value of $\chi$ but it is
$+\infty$ if $\chi\not\in [\chi_*,\chi^*]$ (while $I_{[\chi_*,\chi^*]}(\chi)=0$
if $\chi\in[\chi_*,\chi^*]$). In our approach, instead, we will show that the
constraint on $\chi$ is ensured by the evolution of the system, i.e., it will be
proved that the equations of the system are somehow {\em consistent} as they
yield $\chi_* \leq \chi\leq \chi^*$. The proof of this property of our model
will be detailed in the sequel and relies on a maximum principle argument.

\step Pseudo-potential of dissipation

Secondly, we introduce the pseudo-potential of dissipation $\Phi$ (see
\cite{moreau}) that depends on the dissipative phase variables  $\chi_t$ and
$\nabla\theta$. Let us just comment on the choice of these dissipative
variables: $\chi_t$ is related to microscopic transformations which are responsible
for the phase transition, i.e. for the evolution of the order structure of the
system, while $\nabla\theta$ is concerned with the heat flux. For the sake of
completeness, let us recall that the pseudo-potential of dissipation $\Phi$ is
non-negative, convex with respect to the dissipative variables, and it attains
its minimum $0$ for a null dissipation, that is when
$(\chi_t,\nabla\theta)=(0,{\bf 0})$. We prescribe
\begin{equation}
\Phi(\chi_t,\nabla\theta)=\frac \mu 2|\chi_t|^2+\frac{\lambda}{2\theta}|\nabla\theta|^2
\end{equation}
with $ \mu $ and $\lambda $ denoting positive coefficients. 

\step Constitutive relations

Hence, constitutive relations can be written for $B,\, {\bf H}, \, s,\, {\bf
Q}$. They are recovered from the free energy (for non dissipative contributions)
and the 
pseudo-potential of dissipation (for dissipative parts). We have
\begin{equation}\label{defB}
B=\frac{\partial\Psi}{\partial\chi}+\frac{\partial\Phi}{\partial\chi_t}=\theta_c
F'(\chi)+\theta G'(\chi)+\mu\chi_t
\end{equation}
and
\begin{equation}\label{defH}
{\bf H}=\frac {\partial\Psi}{\partial(\nabla\chi)}=\nu\nabla\chi
\end{equation}
as well as
\begin{equation}\label{defS}
s=-\frac{\partial\Psi}{\partial\theta}=c_0\theta-G(\chi)
\end{equation}
and
\begin{equation}\label{defQ}
{\bf Q}=-\frac{\partial\Phi}{\partial(\nabla\theta)}=-\frac{\lambda}{\theta}
\nabla \theta=- \lambda\nabla\log\theta.
\end{equation}
Let us point out that the choice of the free energy \eqref{free-en} leads to a
linear contribution for the temperature in \eqref{defS}. This will yield
sufficient regularity on the solution, from which we will be able to prove
uniqueness. We also point out that the term $-(c_0/2) \theta^2$ could be seen
as a first order approximation of the following, well-known, form of the energy
potential 
$$
\Psi (\theta, {}\cdots) =-c_0\theta\log\theta + {}\ldots
$$
used, e.g., in \cite{bcf}. In this case, the entropy $s$ would be related to the
temperature $\theta$ through a logarithmic nonlinearity. Note that, on the
contrary, a logarithmic nonlinearity forcing $\theta$ to be strictly positive is
present in our expression \eqref{defQ} for ${\bf Q}$.

\step From energy balance to our equation

The reader may be curious about the derivation of \eqref{entropy} from
\eqref{energybal}.
Let us recall constitutive relations \accorpa{defB}{defQ} and the well known
Helmoltz relation
\begin{equation}\label{defEelena}
 e=\Psi+\theta s .
\end{equation}
By applying the chain rule in \eqref{energybal}, some terms cancel and  
one can rewrite the energy balance~as
\begin{equation}\label{nuovaener}
\theta ( s_t+\div {\bf Q}-R) = \left(\frac{\partial\Phi}{\partial\chi_t},
\frac{\partial\Phi}{\partial(\nabla
\theta)}\right)\cdot(\chi_t,\nabla\theta)
\end{equation}
letting $ {\bf Q}={\bf q}/ \theta $ and  $R= r / \theta $.
We point out that the fact that $\Phi$ is convex, l.s.c., proper, non-negative
and it attains its minimum $0$ when $(\chi_t,\nabla\theta)=(0,{\bf 0})$ ensure
that the right hand side of \eqref{nuovaener} is non-negative. As $\theta>0$,
the Clausius-Duhem inequality
\begin{equation}\label{CD}
 s_t+\div {\bf Q}-R\geq 0.
\end{equation}
complies with \eqref{nuovaener}. 
Hence, dividing \eqref{nuovaener} by $\theta$ and neglecting the resulting
contribution on the right hand side, equation \eqref{entropy} follows.

\step Nonlinearity in the source term

We stress that the entropy source $R$ in \eqref{entropy} is related to the
source $r$ appearing in \eqref{energybal} by $R= r / \theta$. Thus, 
it seems reasonable to include  in our analysis the possibility 
for $R$ to depend on the temperature (and
possibly to present some singularities). This is one of the features of our paper: 
we actually deal with entropy sources (positive or negative, as sinks) 
depending also on the temperature $\theta$. Indeed, possible choices 
satisfying our assumptions are (cf. the later Remark~\ref{Exbetaepi})
$$
 R(x,t,\theta)=\frac{R_1(x,t)}{\theta^2}-R_2(x,t)
$$
which would correspond to $ r(x,t,\theta)= (R_1(x,t)/{\theta}) - R_2(x,t)\theta$, or 
$$
R(x,t,\theta)=R_3(x,t)\theta-R_4(x,t)
$$
which can be viewed as a linearization of $R$ around some equilibrium value of $\theta$. 
In such cases, the possible data $R_1, \, R_2 $ or $R_3, \, R_4 $  are smooth enough and at 
least $R_1$ should be non-negative throughout~$\Omega\times (0,T)$.

\step System of PDEs and initial-boundary value problem 

Now, combining constitutive relations \accorpa{defB}{defQ} with \eqref{eqBH} and
\eqref{entropy} leads to the following PDE system
\begin{align}\label{eqI}
&c_0\theta_t-G'(\chi)\chi_t-\lambda\Delta\log\theta=R(x,t,\theta)\\
\label{eqII}
&\mu\chi_t-\nu\Delta\chi +F'(\chi)\theta_c+G'(\chi)\theta=0
\end{align}
which is addressed in $Q:=\Omega\times(0,T)$.
Then, concerning boundary conditions, prescribed in $\Gamma\times(0,T)$,
we fix a Dirichlet condition for the temperature and a Neumann homogeneous
condition for the phase parameter
\begin{equation}\label{condbord}
\log\theta=\log\theta_\Gamma,\quad\partial_n\chi=0.
\end{equation}
Finally, initial conditions are set in $\Omega$
\begin{equation}\label{condin}
\theta(0)=\theta_0,\quad\chi(0)=\chi_0.
\end{equation}
For the sake of simplicity, in the mathematical analysis performed in
subsequent sections we will take the physical constants $c_0, \, \lambda, \, 
\mu, \, \nu , \, \theta_c $ all equal to 1.   

Let us point out in our model is that we are
dealing with a non-smooth entropy source $R(\theta)$ 
in \eqref{eqI} that is assumed to be
increasing with respect to $\theta \in (0, +\infty)$ up to Lipschitz
perturbations. This fact turns out to be interesting both for analytical and
modelling aspects. Indeed, we are able to treat a diffusive equation for the
temperature with nonlinear and singular diffusion 
along with a non-smooth contribution in source term. Then, for modelling
aspects, we can figure the entropy source $R$ directly as $ {r(\theta)}/\theta$
(compare \eqref{nuovaener} with \eqref{energybal}). The second aspect we notice is
that the physical property $\chi\in[\chi_*,\chi^*]$ comes as a consequence of our
analysis and then it is ensured by the evolution of the system itself and not
required a priori. 

Under suitable assumptions on the data and on the regularity of the involved
nonlinearities, we can prove existence of a solution in any time interval
$(0,T)$ for a weak formulation of our system \accorpa{eqI}{condin}. 
Moreover, in a wide setting we are able to show that
also the solution component $\theta $ is bounded by adapting a Moser type
argument. Finally, if $R(x,t,\theta)$ is uniformly Lipschitz continuous with
respect to $\theta$, we prove uniqueness of the solution. 

\step Phase field model with special heat flux law

What is also interesting in our contribution is that system \accorpa{eqI}{eqII}
may be interpreted as a nonlinear Caginalp phase-field model (see~\cite{cag}) with special
heat flux law. Indeed, take the following energy functional 
\begin{equation}\label{defFPF}
{\Upsilon}(\theta,\chi,\nabla\chi)=-\frac 1{2\theta}+\theta
F(\chi)+G(\chi)+\frac{\nu\theta}2|\nabla\chi|^2.
\end{equation}
Hence, arguing as in the derivation of the Penrose-Fife model (cf. \cite{Pf, BS,
fipe}), we recover the entropy $s$, the internal energy $e$, and the phase field
equation (for $\chi$) as
\begin{align}\label{defSPF}
&s=-\frac{\partial{\Upsilon}}{\partial\theta}=-\frac1{2\theta^2}-F(\chi)-\frac
\nu 2|\nabla\chi|^2\\ \label{defEPF}
&e={\Upsilon}+\theta s=-\frac 1\theta+G(\chi)\\ \label{eqchiPF}
&\mu\chi_t+\frac{\delta H}{\delta \chi}=0,\quad H=\int_\Omega \frac
1\theta{\Upsilon}
\end{align}
where ${\delta H}/{\delta \chi}$ denotes the variational 
(Fr\'echet) derivative of the functional $H.$
Then, we can write the first principle neglecting microscopic motions and forces
(here $\widetilde r$ represents a mere source term)
\begin{equation}\label{IprinC}
e_t+\div{\bf q}=\widetilde r .
\end{equation}
We now let 
\Beq \label{hfl}
{\bf q}=-\lambda\nabla\log\theta 
\Eeq 
and  this special choice of the heat flux can be compared with the ones studied
in the papers \cite{cl, cls}, where Penrose-Fife models with special heat flux
laws have been investigated.  Then, combine \eqref{eqchiPF} and \eqref{IprinC}
with \accorpa{defFPF}{defEPF} and 
\eqref{hfl}, so to obtain the following system
\begin{align}\label{nuovosistI}
&u_t-G'(\chi)\chi_t-\lambda \Delta\log u=-\widetilde r\\
\label{nuovosistII}
&\mu\chi_t-\nu\Delta\chi+F'(\chi)+uG'(\chi)=0
\end{align}
in which $u$ has now the physical meaning of $ 1/ \theta$. Note that system
\accorpa{nuovosistI}{nuovosistII} is formally equivalent to \accorpa{eqI}{eqII}
(with $c_0 =\theta_c =1$ and $- \widetilde r$ in place of $R$).
Thus, it turns out that our analysis applies to the Caginalp system
\accorpa{nuovosistI}{nuovosistII} in which non-smooth heat sources depending on
the temperature are admitted. 

\step Plan of the paper

We conclude the Introduction by giving an outline of the paper. 
In Section~\ref{MainResults} we set our problem
in a Sobolev spaces framework and make precise assumptions on the data and the
involved functions, also stating main results. In Section~\ref{Approx} we
introduce an approximating problem and prove  the existence of solutions for it.
Next, Section~\ref{Existence} is devoted to derive some 
a priori estimates and to pass to the limit as the approximating parameter goes
to 0. A maximum principle is also checked to show that the constraint on $\chi$
is always satisfied. Boundedness of $\theta$ and uniqueness proofs are the
subjects of the last two Sections~\ref{Regularity} and 
\ref{Uniqueness}, respectively.


\section{Main results}
\label{MainResults}
\setcounter{equation}{0}

In this section, we carefully describe the problem we are
going to deal with and state our results.
First of all, we introduce the notation regarding the domain in which
the evolution is considered.
In the sequel, $\Omega$~is a bounded open set in~$\erre^3$,
whose boundary $\Gamma$ is assumed to be of class~$C^2$.
Moreover, $\dn$ is~the (say, outward) normal derivative on~$\Gamma$.
Given a finite final time~$T$,
we set for convenience
\Beq
  Q_t := \Omega \times (0,t) \quad \hbox{for every $t\in(0,T]$}
  \aand
  Q := Q_T \,.
  \label{defQt}
\Eeq
Next, we describe the structure of our system.
We are given four constants
$\thetamin,\,\thetamax,\,\chimin,\,\chimax\in\erre$ such~that
\Beq
  0 < \thetamin \leq 1 \leq \thetamax
  \aand
  \chimin < \chimax
  \label{hpcostanti}
\Eeq
and four functions $F$, $G$, $\Beta$, and~$\pi$
\Beq
  F, G : \erre \to \erre, \quad
  \beta : Q \times (0,+\infty) \to \erre,
  \aand
  \pi : Q \times \erre \to \erre
  \non
\Eeq
satisfying
\Bsist
  \leftone F, G \in C^2(\erre), \quad
  \hbox{$F$ is bounded from below and $G$ is nonnegative}
  \label{hpFG}
  \\
  \leftone F', G' \leq 0 \quad \hbox{in $(-\infty,\chimin)$}, \aand
  F', G' \geq 0 \quad \hbox{in $(\chimax,+\infty)$}
  \label{hpFGbis}
  \\
  \leftone \hbox{$\beta$ is \Lip\ continuous in
      $Q\times[\delta,1/\delta]$ for every $\delta\in(0,1)$}
  \label{hpbetalip}
  \\
  \leftone \hbox{$\betax\,$, $\betat\,$, $\beta'$, and $\pi$
      are Carath\'eodory functions},
  \quad \hbox{with the notation}
  \label{hpcarath}
  \\
  \leftone \betax(x,t,r) := \nabla \beta(x,t,r),
  \quad \betat(x,t,r) := \dt\beta(x,t,r),
  \quad \beta'(x,t,r) := \partial_r \beta(x,t,r)
  \qquad
  \label{defbetax}
  \\
  \leftone 0 \leq \beta'(x,t,r) \leq \betau(r)
  \non
  \\ 
  \leftone \quad {}
  \hbox{\aaQ, every $r\in\erre$, and some $\betau\in C^0(0,+\infty)$}
  \label{hpbeta}
  \\
  \leftone |\betax(x,t,r)| + |\betat(x,t,r)|
  \leq \cbeta (1 + |\beta(x,t,r)|)
  \non
  \\
  \leftone \quad {}
  \hbox{\aaQ, every $r\in(0,+\infty)$, and some $\cbeta\in[0,+\infty)$}
  \label{hpbetaxt}
  \\
  \leftone \beta(x,t,1) = 0 \quad \hbox{for every $(x,t)\in Q$}.
  \label{hpbetanulla}
  \\
  \leftone |\pi(x,t,r)| \leq \lambda |r| + \piz(x,t)
  \non
  \\
  \leftone \quad {}
  \hbox{\aaQ, every $r\in\erre$,
    some $\lambda\in [0,+\infty)$, and some $\piz\in\LQ2$}.
  \label{hppi}
\Esist
\def\HPstruttura{\eqref{hpcostanti}--\eqref{hppi}}%
Furthermore, we set for convenience
\Beq
  \Beta(x,t,r) := \int_1^r \beta(x,t,s) \, ds
  \quad \hbox{for every $(x,t,r)\in Q\times(0,+\infty)$} .
  \label{defprimitbeta}
\Eeq
Then, observing that (by \eqref{hpbeta}) $\beta(x,t,\cpto)$ is nondecreasing, 
it turns out that $\Beta$~is nonnegative and convex with respect 
to the third variable.

\Brem
\label{Rembeta}
We note that the second inequality of \eqref{hpbeta} and \eqref{hpbetanulla}
imply~that
\Beq
  |\beta(x,t,r)|
  \leq \betaz(r)
  := \left| \int_1^r \betau(s) \, ds \right|
  \quad  \forall \ (x,t,r) \in Q\times(0,+\infty).
  \label{dahpbeta}
\Eeq
Therefore, we see by \eqref{hpbetaxt} that even $\betax$ and $\betat$
satisfy an analogous inequality and infer that \eqref{hpbetalip} follows
from the other assumptions.
We have written \eqref{hpbetalip} in advance to give a meaning to the
pointwise values of~$\beta$.
\Erem

\Brem
\label{Exbetaepi}
In the applications just the difference $R=\pi-\beta$ is interesting
and we observe that the choice
\Beq
  R(x,t,r)
  = \frac {R_1(x,t)} {r^2} + R_2(x,t)
  \quad \hbox{with $R_1\geq0$}
  \label{exuno}
\Eeq
is included in our assumption.
Indeed, we can take $\beta(x,t,r)=R_1(x,t)-R_1(x,t)/r^2$ for $r>0$
and $\pi(x,t,r)=R_1(x,t)+R_2(x,t)$, so that $R=\pi-\beta$ and
\eqref{hpbetanulla} hold at the same time.
However, we have to assume that $R_1$ is \Lip\ continuous and
that $|\nabla R_1|+|\dt R_1|\leq c R_1$ in~$Q$ for some constant~$c$ in
order to fulfill~\eqref{hpbetaxt}.
As $R_1$ is bounded and nonnegative, \eqref{hpbeta}~holds as well.
Moreover, we ask that $R_2\in\LQ2$ in view of~\eqref{hppi}.
A~different choice for $R$ is the following
\Beq
  R(x,t,r) = R_3(x,t) \, r - R_4(x,t)
  \label{exdue}
\Eeq
and can be obtained in our setting by simply taking $\beta=0$ and~$\pi=R$.
Here, we can assume $R_3\in\LQ\infty$ and $R_4\in\LQ2$.
\Erem

\Bnot
\label{Notaz}
Let $I$ be a real interval and $\psi:Q\times I\to\erre$ be a Carath\'eodory
function.
We use the same symbol $\psi$ to denote the operator acting on measurable
functions on $Q$ as follows.
If $v:Q\to\erre$ is measurable
\Beq
  \psi(v)
  \quad \hbox{denotes the function} \quad
  (x,t) \mapsto \psi(x,t,v(x,t)), \ (x,t) \in Q.
  \label{notazQ}
\Eeq
Note that $\psi(v)$ actually is measurable due
to the Carath\'eodory assumption on~$\psi$.
Similar definitions and symbols are used for functions depending on
the space variable.
Namely, if $v:\Omega\to\erre$ is measurable
\Beq
  \psi(t,v)
  \quad \hbox{denotes the function} \quad
  x \mapsto \psi(x,t,v(x)), \ x \in \Omega
  \label{notazO}
\Eeq
\aat.
We obtain a time dependent operator.
As a consequence, if $v\in\LQ2$, the symbol $\psi(t,v(t))$ denotes the
measurable function $x\mapsto\psi(x,t,v(x,t))$, $x\in\Omega$.
Notation \accorpa{notazQ}{notazO} will be used, in particular,~with (some
of the functions listed below will be introduced later~on)
\Beq
  \psi = \beta,\ \Beta,\ \betaeps,\ \Betaeps,\
  \betax,\ \betat,\ \beta',\ \betaepsx,\ \betaepst,\ \betaeps',\
  \pi.
  \label{elencopsi}
\Eeq
Furthermore, we~set
\Beq
  H := \Ldue, \quad V := \Huno, \quad \Vz := \Hunoz, \quad
  W := \graffe{v \in \Hdue : \ \dn v = 0}.
  \qquad
  \label{defspazi}
\Eeq
We endow $H$, $V$, and $W$
with their usual scalar products and norms and use a self-explaining
notation, like $\normaV\cpto$.
For the sake of simplicity, the same symbol will be
used both for a space and for any power of~it.
It is understood that
$H$ is embedded in $\Vp$ in the usual~way, i.e.,
$\<u,v>=(u,v)$ for every $u\in H$ and~$v\in\Vz$,
where $\<\cpto,\cpto>$ is the duality pairing
between $\Vp$ and~$\Vz$ and $(\cpto,\cpto)$ is the inner product of~$H$.
\Enot

Now, we list our assumptions on the boundary and initial data.
We are given three functions $\thetaG$, $\thetaz$, and $\chiz$ such~that
\Bsist
  \leftone \thetaG \in \L2{\HxG{1/2}} \cap \H1{\HxG{-1/2}},
  \quad \thetamin \leq \thetaG \leq \thetamax \quad \aeS
  \qquad\quad
  \label{hpthetaG}
  \\
  \leftone \thetaz \in H, \quad
  \thetamin \leq \thetaz \leq \thetamax \quad \aeO
  \label{hpthetaz}
  \\
  \leftone \chiz \in V, \quad
  \chimin \leq \chiz \leq \chimax \quad \aeO
  \label{hpchiz}
\Esist
\def\HPdati{\eqref{hpthetaG}--\eqref{hpchiz}}%
where $\thetamin$, $\thetamax$, $\chimin$, and $\chimax$ are introduced
in~\eqref{hpcostanti}.

At this point, we are ready to state our problem.
We look for a triplet $\soluzione$ satisfying the regularity conditions
and the equations listed below.
\Bsist
  && \theta \in \L\infty H , \quad
   \theta > 0  \quad \aeQ, \aand \ln\theta\in\L2V
  \label{regtheta}
  \\
  && \chi \in \L2W \cap \H1H
  \label{regchi}
  \\
  && G(\chi), \, F'(\chi), \, G'(\chi) \in \LQ2
  \label{regfunzchi}
  \\
  && \xi \in \LQ2
  \label{regxi}
  \\
  && \dt (\theta - G(\chi)) \in \L2\Vp
  \label{regdt}
  \\
  && \dt (\theta - G(\chi)) - \Delta\ln\theta + \xi
  = \pi(\theta)
  \quad \hbox{in $\L2\Vp$}
  \aand
  \xi = \beta(\theta)
  \qquad
  \label{prima}
  \\
  && \dt\chi - \Delta\chi + F'(\chi) + G'(\chi) \, \theta = 0
  \quad \aeQ
  \label{seconda}
  \\
  && \ln\theta = \ln\thetaG \quad \aeS
  \label{bc}
  \\
  && \bigl(\theta - G(\chi) \bigr)(0) = \thetaz - G(\chiz)
  \aand
  \chi(0) = \chiz \,.
  \label{cauchy}
\Esist
\def\regsoluz{\eqref{regtheta}--\eqref{regdt}}%
\def\pbl{\eqref{prima}--\eqref{cauchy}}%
Even though $\xi$ is a known function of~$\theta$,
we refer to the triplet $\soluzione$ instead of the pair $(\theta,\chi)$
when we speak of a solution, just for convenience.
Moreover, we note that \eqref{regtheta} and \eqref{regfunzchi} yield
$G'(\chi)\theta\in\L2\Luno$. 
However, we see by comparison in \eqref{seconda} that $G'(\chi)\theta\in\LQ2$.
Next, we observe that $G(\chiz)$ makes sense in~$\Linfty$.
Furthermore, we point out that the first condition in \eqref{cauchy}
reduces to $\theta(0)=\thetaz$ whenever one knows that $G(\chi)\in\C0\Vp$.
Actually, some additional smoothness for $G(\chi)$ as well as for $F'(\chi)$ and~$G'(\chi)$
surely holds if the nonlinearities satisfy some growth conditions,
thanks to~\eqref{regchi}.
The same is trivially true whenever $\chi$ is bounded, 
and our existence result stated below ensures such a property.
Finally, \eqref{regchi}~itself entails the
homogeneous Neumann boundary condition
for~$\chi$ (see~\eqref{defspazi}).
Here is our first result.

\Bthm
\label{Esistenza}
Let \HPstruttura\ and \HPdati\ be fulfilled.
Then, there exists a triplet $\soluzione$ satisfying the regularity
requirements \regsoluz\ and solving problem \pbl.
Moreover, every solution $\soluzione$ fulfils the inequalities
\Beq
  \chimin \leq \chi \leq \chimax \quad \aeQ.
  \label{massimo}
\Eeq
In particular, $\chi$~is bounded.
\Ethm

One can wonder whether the component $\theta$ of the solution
is bounded as~well.
Actually, such a property holds whenever we reinforce
the assumption on the structure of our system a little, namely
\Beq
  \piz \in \LQ q
  \quad \hbox{for some $q>5/2$}.
  \label{hpq}
\Eeq
We can state the following result.

\Bthm
\label{Regolarita}
Assume \eqref{hpq} in addition to the hypotheses of
Theorem~\ref{Esistenza}.
Then, the component $\theta$ of any solution $\soluzione$ to problem \pbl\
is bounded.
\Ethm

\Brem
\label{Piureg}
From Theorem~\ref{Esistenza} it follows that $G(\chi)$, $F'(\chi)$, and $G'(\chi)$
are smoother, namely
\Beq
  G(\chi), F'(\chi), G'(\chi) \in \LQ\infty \cap \L2V \cap \H1H
  \label{piuregfunzchi}
\Eeq
due to \eqref{hpFG} and~\eqref{regchi}.
In particular, we deduce~that
\Beq
  \dt\theta \in \L2\Vp
  \label{regdttheta}
\Eeq
by comparison in~\eqref{regdt}.
On the other hand, using the regularity conditions $\ln\theta\in\L2V$ and
$\theta\in\LQ\infty$ (see~\eqref{regtheta} and Theorem~\ref{Regolarita}),
we see that $\nabla\theta^m=m\theta^m\nabla\ln\theta$ for every $m>0$,
whence
\Beq
  \theta^m \in \L2V
  \quad \hbox{for every $m\in(0,+\infty)$}.
  \label{thetapiureg}
\Eeq
Moreover, as $F'(\chi)$ and $G'(\chi)$ are bounded, from \eqref{regchi} and
\eqref{seconda} we see that $\dt\chi-\Delta\chi$ is bounded too.
Hence, we~have
\Beq
  \chi \in \L p{\Wx{2,p}} \cap \W{1,p}{\Lx p}
  \quad \hbox{for every $p\geq 1$}
  \label{chipiureg}
  \Eeq
thanks to the general theory of linear parabolic equations.
\Erem

Finally, we state our uniqueness result.
As it often happens for doubly nonlinear problems, uniqueness for solutions
cannot be proved unless more restrictive assumption on the structure are
made.
In particular, we cannot allow a singular $\beta$ like~\eqref{exuno}.

\Bthm
\label{Unicita}
Assume that the hypotheses of Theorem~\ref{Esistenza} are fulfilled.
Moreover, assume that $R:=\pi-\beta$ satisfies the uniform Lipschitz
condition
\Beq
  |R(x,t,r) - R(x,t,s)| \leq \lambdap |r-s|
  \quad \hbox{for every $(x,t)\in Q$ and $r,s\in(0,+\infty)$}
  \label{pimenobetalip}
\Eeq
for some $\lambdap\geq0$.
Then, problem \pbl\ has at most one solution $\soluzione$
satisfying $\theta\in\LQ\infty$.
\Ethm

\Bcor
\label{Riassuntorisultati}
If \eqref{hpq} and \eqref{pimenobetalip} are fulfilled at the same time
in addition to the hypotheses of Theorem~\ref{Esistenza},
then problem \pbl\ has exactly one solution $\soluzione$
such that $\theta$~is bounded.
Moreover, \eqref{massimo}, \eqref{thetapiureg}, and \eqref{chipiureg}
are satisfied.
\Ecor

The paper is organized as follows.
The next section deals with approximating problems.
Theorem~\ref{Esistenza} is proved in Section~\futuro\ref{Existence}
and our argument relies on some a~priori estimates on the approximate solutions
and on monotonicity and compactness methods.
Section~\futuro\ref{Regularity} is devoted to prove Theorem~\ref{Regolarita}
and uses a Moser type procedure.
Finally, the proof of Theorem~\ref{Unicita} is given in the last section.

\medskip

In our proofs, we use the \wk\ inequalities
we are going to recall.
As $\Omega\subset\erre^3$ is bounded and smooth,
the Poincar\'e inequality
\Beq
  \normaV v \leq \cO \normaH{\nabla v}
  \quad \hbox{for every $v\in\Vz$}
  \label{poincare}
\Eeq
holds, and the space $V$ is continuously embedded in $\Lx 6$, i.e.,
\Beq
  \norma v_{\Lx 6} \leq \cO \normaV v
  \quad \hbox{for every $v\in V$}.
  \label{sobolev}
\Eeq
For the sake of completeness,
we recall a related embedding result for parabolic spaces
(see, e.g., \cite[formula~(3.2), p.~8]{DiBen}).
For $m\geq1$, we~have
\Beq
  \L\infty{\Lx m} \cap \L2\Vz \subset \LQ{q(m)}
  \quad \hbox{where} \quad
  q(m) := \frac 23 \, (m+3)
  \non
\Eeq
the embedding being continuous, i.e.,
\Bsist
  && \norma v_{\LQ{q(m)}}
  \leq \cOTm \Bigl( \norma v_{\L\infty{\Lx m}} + \norma{\nabla v}_{\LQ2} \Bigr)
  \non
  \\
  && \quad \hbox{for every $v\in\L\infty{\Lx m} \cap \L2\Vz$}.
  \label{disparabolica}
\Esist
In particular, we observe~that
\Bsist
  && \L\infty H \cap \L2\Vz \subset \LQ{10/3}
  \label{immdue}
  \\
  && \L\infty{\Lx 3} \cap \L2\Vz \subset \LQ 4
  \label{immtre}
\Esist
and that the corresponding estimates \eqref{disparabolica} hold.

The above inequalities are widely used in the sequel,
as well as the elementary Young inequality
\Beq
  ab \leq \delta a^p + c_{\delta,p} \, b^{p'}
  \qquad \forall\  a,b\geq 0 \qquad \forall\  \delta>0
  \label{elementare}
\Eeq
where $p,p'>1$ satisfy $(1/p)+(1/p')=1$
and $c_{\delta,p}:=(p')^{-1}(\delta p)^{-p'/p}$.

We conclude this section by stating a general rule
we use for constants,
in order to avoid a boring notation.
Throughout the paper, the symbol $c$ stands for different constants which depend only
on~$\Omega$, on the final time~$T$, and on the constants and the norms of
the functions involved in the assumptions of either our statements
or our approximation.
In particular, $c$~is independent of the approximation parameter~$\eps$
we introduce in the next section.
A~notation like $\cd$ allows the constant to depend on the positive
parameter~$\delta$, in addition.  Hence, the meaning of $c$ and $\cd$ might
change from line to line and even in the same chain of inequalities.  On
the contrary, we use different symbols (see, e.g.,~\eqref{poincare})
to~denote precise constants which we
could refer~to.
By the way, several of such constants
could be the same (like in \eqref{poincare} and in~\eqref{sobolev}),
since sharpness is not needed.


\section{Approximating problems}
\label{Approx}
\setcounter{equation}{0}

This section contains a preliminar work in the direction of proving
Theorem~\ref{Esistenza} and deals with a suitable appoximation of problem
\pbl.
Namely, we replace the strong nonlinearities that appear in
equation~\eqref{prima} with smooth functions depending on the parameter
$\eps\in(0,1)$.
First of all, we see the logarithm $\ln$ as a maximal monotone operator
in~$\erre$ with domain~$(0,+\infty)$.
Precisely, for $r\in\erre$, the set $\ln r$ is the singleton $\{\ln r\}$ if $r>0$, while it
is empty if $r\leq0$.
Then, we can consider the function $\Lneps:\erre\to\erre$ defined
as~follows
\Beq
  \Lneps r := \eps r + \lneps r
  \quad \hbox{where} \quad
  \hbox{$\lneps$ is the Yosida \regulariz ation of $\ln$}.
  \label{defLneps}
\Eeq
We note that $\lneps$ is monotone and \Lip\ continuous with
constant~$1/\eps$ (see, e.g., \cite[Prop.~2.6, p.~28]{Brezis}).
Thus, $\eps\leq\Lneps'(r)\leq\eps+(1/\eps)$ for every $r\in\erre$.
The functions $\lneps$ and $\Lneps$ act on $L^2$-spaces as well and
Notation~\ref{Notaz} is extended to~them.
Moreover, we replace $F$ and $G$ by new functions, still termed $F$
and~$G$, satisfying
\Beq
  \hbox{$F'$ and $G'$ are bounded}
  \label{lipFG}
\Eeq
in addition to~\eqref{hpFG}--\eqref{hpFGbis}.
Indeed, we can arbitrarily modify $F$ and $G$ outside $[\chimin,\chimax]$
due to the last part of Theorem~\ref{Esistenza}.
Finally, we replace the possibly singular function $\beta$ by a $C^\infty$
function~$\betaeps$, in order to justify the chain rules we have to use.
We proceed by extension, truncation, and \regulariz ation.  For $\eps>0$,
we~set
\Beq
  \Omega_\eps
  := \graffe{x\in\Omega: \ \dist(x,\Gamma)<\eps}
  \aand
  \Omega'_\eps
  := \graffe{x\in\erre^3\setminus\overline\Omega: \ \dist(x,\Gamma)<\eps}.
  \non
\Eeq
As $\Omega$ is of class~$C^2$, there exists $\eps_0\in(0,1)$ such that, for
every $x\in\Omega'_{\eps_0}$, there exists a unique point
$\tilde x\in\Omega_{\eps_0}$ satisfying
\Beq
  x' := \frac {x+\tilde x} 2 \in \Gamma
  \aand
  \hbox{$x-\tilde x$ is orthogonal to $\Gamma$ at $x'$}
  \label{riflessione}
\Eeq
the correspondence $x\mapsto\tilde x$ being a bi-\Lip\ diffeomorphism of
class $C^1$ from $\Omega'_{\eps_0}$ onto~$\Omega_{\eps_0}$.
Then, we define $\Omegatilde$ and the extension-by-reflection
operator
$\widetilde{\cpto}:\Linfty\to L^\infty(\Omegatilde)$
as~follows
\Bsist
  && \Omegatilde := \overline\Omega \cup \Omega'_{\eps_0}
  \aand
  \hbox{for $v\in\Linfty$ and a.a.\ $x\in\Omegatilde$ we set}
  \non
  \\
  && \widetilde v(x) = v(x)
  \quad \hbox{if $x\in\Omega$}
  \aand
  \widetilde v(x) = v(\tilde x)
  \quad \hbox{if $x\in\Omega'_{\eps_0}$}.
  \label{extOmega}
\Esist
Next, we extend further by reflection as well.
We define $\Qtilde$ and
$\calE:\LQ\infty\to L^\infty(\Qtilde)$ as follows
\Bsist
  && \Qtilde := \Omegatilde \times (-T,2T)
  \aand
  \hbox{for $v\in\LQ\infty$ and a.a.\ $(x,t)\in\Qtilde$ we set}
  \non
  \\
  && (\calE v)(x,t) = \widetilde v(x,t)
  \quad \hbox{if $t\in(0,T)$}, \quad
  (\calE v)(x,t) = \widetilde v(x,-t)
  \quad \hbox{if $t\in(-T,0)$}
  \non
  \\
  && \aand
  (\calE v)(x,t) = \widetilde v(x,2T-t)
  \quad \hbox{if $t\in(T,2T)$}.
  \label{defE}
\Esist
It is clear that the extension operator $\calE$ is linear and continuous.
More precisely, we~have
\Beq
  \supess_{\Qtilde} \calE v = \supess_Q v
  \aand
  \infess_{\Qtilde} \calE v = \infess_Q v
  \quad \hbox{for every $v\in\LQ\infty$}.
  \label{lincontE}
\Eeq
Moreover, one can check~that for every $v\in\LQ\infty$ we~have
\Bsist
  && \calE v \geq 0 \, \hbox{ a.e. in } \Qtilde 
  \quad \hbox{whenever} \quad
  v \geq 0 \, \hbox{ a.e. in } Q 
  \label{posE}
  \\
  && \norma{\nabla \calE v}_{L^\infty(\Qtilde)}
  \leq M \norma{\nabla v}_{\LQ\infty}
  \quad \hbox{if $\nabla v\in\LQ\infty$}
  \label{regEx}
  \\
  && \norma{\dt \calE v}_{L^\infty(\Qtilde)}
  \leq M \norma{\dt v}_{\LQ\infty}
  \quad \hbox{if $\dt v\in\LQ\infty$}
  \label{regEt}
  \\
  && \lip (\calE v) \leq M \lip v
  \quad \hbox{if $v$ is \Lip\ continuous}
  \label{regElip}
\Esist
for some constant~$M$, where $\lip v$ is the \Lip\ constant of~$v$.
At this point, for $\eps \in (0,1)$ we define the operators
$\betatilde,\betatildeeps:\Qtilde\times\erre\to\erre$
by using the extension operator $\calE$ and a truncation procedure
as~follows
\Bsist
  && \betatilde(x,t,r)
  := \bigl( \calE \beta(\cpto,\cpto,r) \bigr) (x,t)
  \label{defbetatilde}
  \\
  && \betatildeeps(x,t,r)
  := \betatilde(x,t,r_\eps)
  \quad \hbox{where} \quad
  r_\eps := \max \{\eps,\min \{r,1/\eps\}\}
  \label{defbetatildeeps}
\Esist
and observe that $\betatildeeps$ is globally \Lip\ continuous.
Indeed, by recalling \eqref{hpbetalip} and setting for convenience
\Beq
  L_\delta := \lip \beta|_{Q\times[\delta,1/\delta]}
  \quad \hbox{for $\delta\in(0,1)$}
  \label{defLdelta}
\Eeq
we clearly see that $ML_\eps$ is a \Lip\ constant for~$\betatildeeps$.
Moreover, as $\calE$ is linear and \eqref{posE} holds, we infer that both
$\betatilde(x,t,\cpto)$ and $\betatildeeps(x,t,\cpto)$ are nondecreasing
on~$\erre$ for every $(x,t)\in\Qtilde$.
Furthermore, both $\betatilde$ and $\betatildeeps$ vanish at $r=1$.
In particular, their values at every $r\in\erre$ have the sign of $r-1$.

Finally, we are ready to define a $C^\infty$ approximation
$\betaeps:Q\times\erre\to\erre$ of~$\beta$.
We \regulariz e $\betatildeeps$ by convolution and restrict the
\regulariz ation we obtain to~$Q\times\erre$.
Namely, we fix a nonnegative $\zeta\in C^\infty(\erre^5)$ supported in the
unit ball $B$ of~$\erre^5$ and \normaliz ed in $L^1(\erre^5)$.
Then, by assuming $\eps_0\leq T$ and $\eps\in(0,\eps_0)$
(such restrictions are not stressed in the sequel, but it is understood
that they are satisfied),
we recall \eqref{defLdelta} and~set
\Bsist
  &&\deps := \frac \eps {1 + L_\eps}
  \aand
  \zetaeps(x,t,r)
  := \deps^{-5} \, \zeta\bigl( (x,t,r)/\deps \bigr)
  \quad \hbox{for $(x,t,r)\in\erre^5$}
  \label{defdeps}
  \\
  && \betaeps(x,t,r)
  := (\betatildeeps * \zetaeps) (x,t,r)
  = \int_{\deps B}
    \betatildeeps(x-y, t-\tau,r-s)
    \, \zetaeps(y,\tau,s)
    \, dy \, d\tau \, ds
  \non
  \\
  && = \int_B
    \betatildeeps(x-\deps y, t-\deps\tau,r-\deps s)
    \, \zeta(y,\tau,s)
    \, dy \, d\tau \, ds
  \quad \hbox{for $(x,t,r)\in Q\times\erre$}.
  \label{defbetaeps}
\Esist
The reason of the above choice of $\deps$ is that we would like to~have
\Beq
  |\betaeps(x,t,r) - \betatildeeps(x,t,r)|
  \leq M \eps
  \quad \hbox{for every $(x,t,r)\in Q\times\erre$}
  \label{erroreregol}
\Eeq
and some constant~$M$.
Actually, \eqref{erroreregol}~holds with the constant $M$ that makes
\eqref{regElip} true, as we show at once.
We have~indeed
\Bsist
  \leftone |\betaeps(x,t,r) - \betatildeeps(x,t,r)|
  = \left| \int_B
    \bigl(
      \betatildeeps(x-\deps y, t-\deps\tau,r-\deps s)
      - \betatildeeps(x,t,r)
    \bigr)
    \, \zeta(y,\tau,s)
    \, dy \, d\tau \, ds
    \right|
  \non
  \\
  \leftone \leq \int_B
    M L_\eps |(\deps y,\deps\tau,\deps s)| \, \zeta(y,\tau,s)
    \, dy \, d\tau \, ds
    \leq M L_\eps \deps
    \leq M \eps
  \non
\Esist
since $ML_\eps$ is a \Lip\ constant for~$\betatildeeps$, as just observed.
With a similar argument, we see~that
\Beq
  \hbox{$\betaeps$ is \Lip\ continuous with}
  \quad \lip \betaeps \leq M L_\eps
  \label{betaepslip}
\Eeq
since such a property holds for~$\betatildeeps$.
In the sequel we use the following more precise~facts
\Bsist
  && \sup_{Q \times [\delta,1/\delta]} |\betaeps|
  \leq \sup_{\Qtilde \times [\delta/2,1/\delta+\delta/2]} |\betatilde|
  \aand
  \lip \betaeps |_{Q \times [\delta,1/\delta]}
  \leq \lip \betatilde |_{\Qtilde \times [\delta/2,1/\delta+\delta/2]}
  \non
  \\
  && \quad \hbox{for} \quad \delta \in (0,1)
  \aand
  \eps < \delta/2 \,.
  \label{lipdelta}
\Esist
Indeed, we have $\deps\leq\eps\leq\delta/2$.
Hence, if $(x,t)\in Q$ and $\delta\leq r\leq1/\delta$,
the values of $\betatildeeps$ in \eqref{defbetaeps}
actually are values of $\betatilde$ at points of the set
$\Qtilde\times[\delta/2,1/\delta+\delta/2]$,
where $\betatilde$ is bounded and \Lip\ continuous.
Therefore, both the supremum and the \Lip\ constant are preserved by the
convolution, since $\zeta$ is \normaliz ed, and \eqref{lipdelta} follow.
Finally, we point out~that
\Beq
  \hbox{$\betaeps(x,t,\cpto)$ is nondecreasing on $\erre$
    for every $(x,t)\in Q$}
  \label{betaepsmon}
\Eeq
since such a property holds for $\betatildeeps$ and $\zeta$ is nonnegative.
Moreover, we set for convenience
\Beq
  \Betaeps(x,t,r) := \int_1^r \betaeps(x,t,s) \, ds
  \quad \hbox{\aaQ\ and every $r\in\erre$}.
  \label{defBetaeps}
\Eeq
Clearly, $\Betaeps$~is convex with respect to the third variable.
Furthermore, as $\betatildeeps(x,t,r)$ and $r-1$ have the same sign
as just observed, we see that \eqref{erroreregol}~implies
\Beq
  \Betaeps(x,t,r)
  \geq \int_1^r \bigl( \betaeps(x,t,s) - \betatildeeps(x,t,s) \bigr) \, ds
  \geq - M \eps |r-1|
  \quad \hbox{for every $(x,t,r)\in Q\times\erre$}.
  \label{belowBetaeps}
\Eeq
Finally, it is clear that the first of \eqref{lipdelta} implies the
analogue for~$\Betaeps$, namely
\Beq
  \sup_{Q \times [\delta,1/\delta]} |\Betaeps|
  \leq \cd
  \quad \hbox{for every $\delta\in(0,1)$}
  \label{bddBetaeps}
\Eeq
for some constant~$\cd$.

At this point, we are ready to state the approximating problem, which
consists in finding a triplet $\soluz\eps$ having the proper regularity and
satisfying
\Bsist
  && \dt (\thetaeps - G(\chieps)) - \Delta\Lneps\thetaeps + \xieps
  = \pi(\thetaeps)
  \quad \hbox{in $\L2\Vp$}
  \aand
  \xieps = \betaeps(\thetaeps)
  \qquad
  \label{primaeps}
  \\
  && \dt\chieps - \Delta\chieps + F'(\chieps) + G'(\chieps) \, \thetaeps = 0
  \quad \aeQ
  \label{secondaeps}
  \\
  && \thetaeps = \thetaG \quad \aeS
  \label{bceps}
  \\
  && \thetaeps(0) = \thetaz
  \aand
  \chieps(0) = \chiz \,.
  \label{cauchyeps}
\Esist
\def\pbleps{\eqref{primaeps}--\eqref{cauchyeps}}%
The following result holds.

\Bthm
\label{Esistenzaeps}
Let the assumptions of Theorem~\ref{Esistenza} be fulfilled.
Moreover, assume \eqref{defLneps}, \eqref{lipFG}, and~\eqref{defbetaeps}.
Then, problem \pbleps\ has a solution satisfying
\Bsist
  && \thetaeps \in \L2V \cap \H1\Vp
  \label{regthetaeps}
  \\
  && \chieps \in \L2W \cap \H1H
  \label{regchieps}
\Esist
\Ethm

We avoid proving Theorem~\ref{Esistenzaeps}
in order not to make the paper too long.
Indeed, the a~priori estimates in Section~\futuro\ref{Existence}
suggest how to proceed.
Anyway, a~rigorous proof could be done by \regulariz ing $\pi$
and using, e.g., a Galerkin procedure.

On the contrary, we prove some auxiliary results regarding
the approximating nonlinearities.
The first formula we state is showed in \cite[Lemma~6.1]{BCFG1}.
We repeat the short proof here, for convenience, and include it in the
following proposition.

\Bprop
\label{Stimelneps}
We have
\Bsist
  && \lneps^{-1}(s) = e^s + \eps s
  \quad \hbox{for every $s\in\erre$}
  \label{invlneps}
  \\
  && \frac {\ln r}{1+\eps} \leq \lneps r \leq \ln r
  \quad \hbox{for every $r\geq1$}
  \label{stimalneps}
  \\
  && \Lneps'(r) \geq 1
  \quad \hbox{for every $r\leq1$}
  \aand
  \Lneps'(r) \geq \frac 1 {2r}
  \quad \hbox{for every $r\geq1$}
  \label{coerclneps}
  \\
  && \lmin \leq \lneps r \leq \lmax
  \quad \hbox{for every $r\in\Istar$}
  \label{stimastar}
\Esist
where we have set $\lmin:=\min\graffe{0,\ln\thetamin}$ and
$\lmax:=\max\graffe{0,\ln\thetamax}$.
Moreover, we~have
\Beq
  \frac 1 {2\thetamax}
  \leq \lneps'(r)
  \leq \frac 2 \thetamin
  \quad \hbox{for every $r\in\Istar$}
  \label{stimastarprimo}
\Eeq
for $\eps$ small enough.
\Eprop

\Bdim
For $r\in\erre$, let $\rhoeps(r)>0$ be defined by the equation
\Beq
  \rhoeps(r) + \eps \ln\rhoeps(r) = r .
  \label{defrhoeps}
\Eeq
Take now any $s\in\erre$. Then, $r:=e^s+\eps s$ satisfies $\rhoeps(r)=e^s$.
On the other hand, we~have
\Beq
  \lneps r = \frac {r - \rhoeps(r)} \eps
  \label{yosidaln}
\Eeq
by definition of Yosida \regulariz ation.
We deduce that $\lneps r=s$ and \eqref{invlneps} follows.
To prove~\eqref{stimalneps}, we observe that for $s\geq0$ we have
$e^{\eps s} \geq 1+\eps s$ and $e^s\geq1$.
We infer~that
\Beq
  e^{(1+\eps)s}
  \geq e^s + \eps s
  \geq e^s .
  \non
\Eeq
If $r\geq 1$, applying this to $s:=\lneps r$
(which is nonnegative since $\lneps 1=0$),
we~obtain
\Beq
  e^{(1+\eps)\lneps r}
  \geq e^{\lneps r} + \eps \lneps r
  \geq e^{\lneps r}
\non
\Eeq
and \eqref{stimalneps} follows from~\eqref{invlneps}
by applying the logarithm.
Next, we prove~\eqref{coerclneps}.
To this aim, we compute $\lneps'(r)$ from
\eqref{yosidaln} and~\eqref{defrhoeps}.
We~have
\Beq
  \lneps'(r) = \frac {1 - \rhoeps'(r)} \eps
  = \frac 1\eps \tonde{ 1 - \frac {\rhoeps(r)} {\rhoeps(r)+\eps} }
  = \frac 1 {\rhoeps(r)+\eps} \,.
  \label{lnepsprimo}
\Eeq
On the other hand,
we observe that \eqref{defrhoeps} and $\rhoeps(r)>1$ imply
$r>\rhoeps(r)\geq1$.
Hence, $\rhoeps(r)\leq1$ for $r\leq1$.
We conclude~that
\Beq
  \Lneps'(r) = \eps + \lneps'(r)
  \geq \eps + \frac 1 {1+\eps}
  = \frac {1+\eps+\eps^2} {1+\eps}
  \geq 1
  \non
\Eeq
for every $r\leq1$.
Assume now $r\geq1$.
Then, $r+\eps\ln r\geq r$, whence $\rhoeps(r)\leq r$.
Accounting for~\eqref{lnepsprimo}, we infer~that
\Beq
  \lneps'(r) \geq \frac 1 {r+\eps} \geq \frac 1 {2r}
  \non
\Eeq
and the second of~\eqref{coerclneps} follows.
To prove~\eqref{stimastar}, we observe that
$\exp\lmin+\eps\lmin\leq\thetamin$ and $\exp\lmax+\eps\lmax\geq\thetamax$.
We deduce that $\exp\lmin+\eps\lmin\leq r\leq\exp\lmax+\eps\lmax$ for every
$r\in\Istar$, and \eqref{stimastar} follows from~\eqref{invlneps}.
Finally, we prove~\eqref{stimastarprimo}.
Owing to~\eqref{lnepsprimo} and to the monotonicity of~$\rhoeps$, we see
that $r\in\Istar$ implies~that
\Beq
  \frac 1 {\rhoeps(\thetamax) + \eps}
  \leq \lneps'(r)
  \leq \frac 1 {\rhoeps(\thetamin) + \eps} \,.
  \non
\Eeq
On the other hand, it is clear that $\rhoeps(r')$ tends to $r'$ as $\eps\to0$ for every
$r'>0$ (see also, e.g., \cite[Prop.~2.6, p.~28]{Brezis}).
Then, \eqref{stimastarprimo}~immediately follows if $\eps$ is small enough.
\Edim

Next, we deal with the analogue of \eqref{hpbetaxt} for~$\betaeps$.
At the same time, we prove an inequality involving~$\Betaeps$.
A~notation like \eqref{defbetax} is extended to such functions.

\Bprop
\label{Propbetaxt}
We have
\Beq
  |\betaepsx(x,t,r)| + |\betaepst(x,t,r)|
  \leq c \bigl( 1 + |\betaeps(x,t,r)| \bigr)
  \aand
  |\Betaepst(x,t,s)|
  \leq c \bigl( \Betaeps(r) + |r| + 1 \bigr)
  \label{disugbetaepsxt}
\Eeq
for every $(x,t,r)\in Q\times\erre$,
some constant~$c$, and $\eps$ small enough.
\Eprop

\Bdim
As far as the first inequality is concerned,
we deal, e.g., with $\betaepst$, since the argument for the space
derivatives is quite similar.
We first prove~that
\Beq
|\dt\betatildeeps(x,t,r)|
\leq M \cbeta \bigl( 1 + |\betatildeeps(x,t,r)| \bigr)
\label{perbetaepst}
\Eeq
for every $(x,t)\in\Qtilde$ and $r\in\erre$, where $M$ and $\cbeta$ are the constants
satisfying \accorpa{regEx}{regElip} and \eqref{hpbetaxt}, respectively.
Assume first $r\in[\eps,1/\eps]$.
Then, \eqref{perbetaepst}~coincides with the analogue for~$\betatilde$
due to~\eqref{defbetatildeeps},
and this easily follows from
\eqref{regEt}, \eqref{defbetatilde}, and~\eqref{hpbetaxt}.
Assume now $r<\eps$.
Then, we~have
\Beq
  |\dt\betatildeeps(x,t,r)|
  = |\dt\betatilde(x,t,\eps)|
  \leq M \cbeta \bigl( 1 + |\betatilde(x,t,\eps)| \bigr)
  = M \cbeta \bigl( 1 + |\betatildeeps(x,t,r)| \bigr)
  \non
\Eeq
just using the above with $r=\eps$.
As the argument is similar if $r>1/\eps$, \eqref{perbetaepst}~is
established.
In order to prove the first of~\eqref{disugbetaepsxt}, we notice~that
\Beq
  \pm \betaepst
  = \pm \dt (\betatildeeps * \zetaeps)
  = \pm (\dt\betatildeeps) * \zetaeps
  \leq c (1 + |\betatildeeps|) * \zetaeps
  = c + c |\betatildeeps| * \zetaeps
  \non
\Eeq
with $c:=M\cbeta$, since \eqref{perbetaepst} holds and the convolution with
the nonnegative \normaliz ed kernel $\zetaeps$ preserves order and
constants.
Therefore, we have to bound the last convolution with the \rhs\ of the
inequality we want to prove.
Assume first $(x,t)\in Q$ and $r\geq1+\eps$.
Then, $\betatildeeps(y,\tau,s)\geq0$ for $(y,\tau)\in\Qtilde$ and
$|s-r|< \deps$ (since $\deps\leq\eps$), and we~have
\Beq
  \bigl( |\betatildeeps| * \zetaeps \bigr) (x,t,r)
  = \bigl( \betatildeeps * \zetaeps \bigr) (x,t,r)
  = \betaeps (x,t,r)
  = |\betaeps (x,t,r)| .
  \non
\Eeq
If $r<1-\eps$ the argument is similar.
Finally, if $|r-1|\leq\eps$, by assuming $\eps\leq1/4$, we have
$\deps\leq1/4$, $r-\deps\geq\max\graffe{\eps,1/2}$, and
$r+\deps\leq\min\graffe{1/\eps,3/2}$, whence
\Beq
  \bigl| \bigl( |\betatildeeps| * \zetaeps \bigr) (x,t,r) \bigr|
  = \bigl| \bigl( |\betatilde| * \zetaeps \bigr) (x,t,r) \bigr|
  \leq \sup_{\Qtilde\times[r-\deps,r+\deps]} |\betatilde|
  \leq \sup_{\Qtilde\times[1/2,3/2]} |\betatilde|
  \non
\Eeq
since $\zetaeps$ is \normaliz ed in~$L^1(\erre^5)$,
and the first of \eqref{disugbetaepsxt} clearly follows.
Then, we easily derive the second inequality.
As $\betaeps$ is monotone with respect to the third variable, we~have
\Bsist
  && |\Betaepst(x,t,r)|
  = \left| \int_1^r \betaepst(x,t,s) \, ds \right|
  \leq c \left| \int_1^r \bigr( 1 + |\betaeps(x,t,s)| \bigr) \, ds \right|
  \non
  \\
  && \leq c |r-1|
  + c \int_1^r \bigl( \betaeps(x,t,s) - \betaeps(x,t,1) \bigr) \, ds
  + c \left| \int_1^r |\betaeps(x,t,1)| \, ds \right|
  \non
  \\
  && \leq c \bigl( \Betaeps(x,t,r) + |r| + 1 \bigr)
  \non
\Esist
since $\betaeps(x,t,1)$ is bounded due to \eqref{lipdelta} with
$\delta=1/2$.
\Edim

The last preliminary remarks we make regard relationships between the
approximating nonlinearities and the boundary datum~$\thetaG$.
In the next section, we need to consider a known smooth function that
coincides with $\thetaG$ on the boundary.
Thus, a~natural choice is the harmonic extension $\thetaH$ of~$\thetaG$.
Precisely, we define $\thetaH:Q\to\erre$ by the conditions
\Beq
  \thetaH(t) \in V, \quad
  \thetaH(t)|_\Gamma = \thetaG(t) ,
  \aand
  \Delta\thetaH(t) = 0
  \quad \hbox{in $\Omega$},
  \quad \aat.
  \label{defthetaH}
\Eeq
Then, assumption \eqref{hpthetaG} and the general theory of harmonic
functions (in~particular, the maximum principle) ensure that
$\thetaH\in\L2V\cap\H1H$ and that the following estimates~hold
\Bsist
  && \norma\thetaH_{\L2V\cap\H1H}
  \leq c \norma\thetaG_{\L2{\HxG{1/2}}\cap\H1{\HxG{-1/2}}}
  \non
  \\
  && \aand
  \thetamin \leq \thetaH \leq \thetamax
  \quad \aeQ.
  \label{regthetaH}
\Esist

\Bprop
\label{StimethetaH}
We have
\Bsist
  && \norma{\Lneps\thetaH}_{\LQ\infty\cap\L2V}
  \leq c
  \label{stimaLnthetaH}
  \\
  && \norma{\betaeps(\thetaH)}_
    {\LQ\infty \cap \L2V \cap \H1H}
  \leq c
  \label{stimabetaH}
\Esist
for some constant $c$ and $\eps$ small enough.
\Eprop

\Bdim
Estimate \eqref{regthetaH} and \eqref{stimastar} imply that
$\lmin\leq\lneps\thetaH\leq\lmax$ \aeQ,
where the notation of Proposition~\ref{Stimelneps} has been used.
Owing to \eqref{stimastarprimo} as well, \eqref{stimaLnthetaH}~immediately
follows.
To prove the $L^\infty$-bound of~\eqref{stimabetaH}, it suffices to recall
\eqref{hpcostanti}, \eqref{dahpbeta}, and the first of~\eqref{lipdelta}.
Finally, we prove the estimate regarding the time derivative
(the~argument for space derivatives is similar).
We~have
\Bsist
  && \norma{\dt\betaeps(\thetaH)}_{\LQ2}
  \leq \norma{\betaepst(\thetaH)}_{\LQ2}
  + \norma{\betaeps'(\thetaH) \, \dt\thetaH}_{\LQ2}
  \non
  \\
  && \leq c \norma{\betaepst(\thetaH)}_{\LQ\infty}
  + \sup_{Q\times[\thetamin,\thetamax]} |\betaeps'| \
    \norma{\dt\thetaH}_{\LQ2}
  \leq c
  \non
\Esist
due to \eqref{disugbetaepsxt}, the $L^\infty$-estimate just proved,
the second of \eqref{lipdelta}, and~\eqref{regthetaH},
provided that $\eps$ is small enough.
\Edim


\section{Existence}
\label{Existence}
\setcounter{equation}{0}

In this section, we prove Theorem~\ref{Esistenza}.
More precisely, we first check the last part of the statement, i.e.,
we show that every solution to problem \pbl\ satisfies the
bounds~\eqref{massimo}.
Then, we consider the approximating problem \pbleps\ taking
the further assumption \eqref{lipFG} into account, and perform a number of a~priori
estimates on its solution~$\soluz\eps$.
Finally, we let $\eps$ tend to zero by using monotonicity and compactness
methods.
This leads to the proof of Theorem~\ref{Esistenza} under the additional
condition~\eqref{lipFG}.
However, it is clear that such a procedure proves Theorem~\ref{Esistenza}
in the general case.
Indeed, the triplet $\soluzione$ we find solves anyone of the problems
obtained by replacing the original functions $F$ and $G$ by other ones that
coincides with the given $F$ and $G$ on~$[\chimin,\chimax]$, due to~\eqref{massimo}.

Hence, we assume~\eqref{lipFG} throghout the present section.
As such a condition is added to our assumptions, we allow the values of the
different constants $c$ to depend even on the \Lip\ constants of $F$
and~$G$, following the general rule explained at the end of
Section~\ref{MainResults}.
Moreover, it is understood that $\eps$ belongs to $(0,1)$ and is even
smaller as in the previous section (see,~e.g.,
Propositions~\ref{Propbetaxt} and~\ref{StimethetaH}).
Finally, $\delta$~is a positive parameter, say $\delta\in(0,1)$, whose
value is chosen according to our convenience.

\step
The maximum principle argument

We prove that every solution
$\soluzione$ to problem \pbl\
satisfies~\eqref{massimo}.
Fix a \Lip\ continuous function $H:\erre\to\erre$ of class $C^1$ such~that
\Beq
  H(r) = 0
  \quad \hbox{if $r\in[\chimin,\chimax]$}
  \aand
  H'(r) > 0
  \quad \hbox{if $r\not\in[\chimin,\chimax]$}.
  \label{hpH}
\Eeq
As $\chi\in\L2V$, it turns out
that $v:=H(\chi)$ is an admissible test function
for~\eqref{seconda}.
Then, we multiply \eqref{seconda} by $H(\chi)$ and integrate over~$Q_t$,
where $t\in(0,T]$ is arbitray.
After integrating by parts, we~obtain
\Beq
  \intQt \dt\chi \, H(\chi)
  + \intQt \nabla \chi \cdot \nabla H(\chi)
  + \intQt \bigl( F'(\chi) + G'(\chi) \, \theta \bigr) \, H(\chi)
  = 0 .
  \label{permassimo}
\Eeq
If $\hat H$ denotes the primitive of $H$ vanishing at~$\chimin$,
we~have
\Beq
  \intQt \dt\chi \, H(\chi)
  = \iO \hat H(\chi(t))
  - \iO \hat H(\chiz)
  = \iO \hat H(\chi(t))
  \geq 0
  \non
\Eeq
due to \eqref{hpchiz} and~\eqref{hpH}.
The next integral of \eqref{permassimo} is nonnegative, obviously, and the
last one has the same property.
Indeed, its integrand is nonnegative since $\theta>0$ and the three
functions $F'$, $G'$, and $H$ have the same sign, due to \eqref{hpFGbis}
and~\eqref{hpH}.
Therefore, all the integrals of \eqref{permassimo} vanish identically.
In particular, we deduce that $\hat H(\chi(t))=0$ \aeO, for every
$t\in[0,T]$, and this clearly implies~\eqref{massimo}.

Now, we start proving estimates on the solution $\soluz\eps$ to problem
\pbleps, as mentioned at the beginning of the present section.
When we shortly say that we test an equation by some function~$v$, we mean
that we test such an equation at time $s$ by~$v(s)$, we integrate first
over~$\Omega$ and then over~$(0,t)$ with respect to~$s$, where $t\in(0,T)$
is arbitrary, and we integrate by parts, if necessary.

\step First a priori estimate

We test \eqref{primaeps} by
$v:=\thetaeps-\thetaH+\delta(\Lneps\thetaeps-\Lneps\thetaH)$.
At the same time, we test \eqref{secondaeps} by~$\dt\chieps$.
Then, we sum the equalities we get to each other and note that two terms
involving $G$ cancel.
After rearranging a little and adding the same integral
to both sides for convenience, we~obtain
\Bsist
  && \frac 12 \iO |\thetaeps(t) - \thetaH(t)|^2
  + \intQt \dt G(\chieps) \, \thetaH
  + \intQt \nabla\Lneps\thetaeps \cdot \nabla\thetaeps
  \non
  \\
  && \quad {}
  + \intQt (\betaeps(\thetaeps) - \betaeps(\thetaH)) (\thetaeps - \thetaH)
  \non
  \\
  && \quad {}
  + \delta \intQt \dt\thetaeps \, \Lneps\thetaeps
  + \delta \intQt |\nabla (\Lneps\thetaeps - \Lneps\thetaH)|^2
  \non
  \\
  && \quad {}
  + \delta \intQt
      (\betaeps(\thetaeps) - \betaeps(\thetaH))
      (\Lneps\thetaeps - \Lneps\thetaH)
  \non
  \\
  && \quad {}
  + \intQt |\dt\chieps|^2
  + \frac 12 \iO |\nabla\chieps(t)|^2
  + \frac 12 \iO |\chieps(t)|^2
  \non
  \\
  && = \frac 12 \iO |\thetaz - \thetaH(0)|^2
  - \intQt \dt\thetaH (\thetaeps - \thetaH)
  + \intQt \nabla\Lneps\thetaeps \cdot \nabla\thetaH
  \non
  \\
  && \quad {}
  - \intQt \betaeps(\thetaH) (\thetaeps - \thetaH)
  + \intQt \pi(\thetaeps) (\thetaeps - \thetaH)
  \non
  \\
  && \quad {}
  + \delta \intQt \dt\thetaeps \, \Lneps\thetaH
  + \delta \intQt G'(\chieps) \, \dt\chieps \,
      (\Lneps\thetaeps - \Lneps\thetaH)
  \non
  \\
  && \quad {}
  - \delta \intQt
      \nabla\Lneps\thetaH
      \cdot \nabla (\Lneps\thetaeps - \Lneps\thetaH)
  - \delta \intQt \betaeps(\thetaH) (\Lneps\thetaeps - \Lneps\thetaH)
  \non
  \\
  && \quad {}
  + \frac 12 \iO |\nabla\chiz|^2
  + \delta \intQt \pi(\thetaeps) (\Lneps\thetaeps - \Lneps\thetaH)
  - \intQt F(\chieps) \dt\chieps
  + \frac 12 \iO |\chieps(t)|^2 .
  \qquad \quad
  \label{perprimastima}
\Esist
Now, we observe with the help of monotonicity
that all terms on the \lhs\ are nonnegative but two of them.
We deal with the first one that needs some treatment.
We~have
\Bsist
  \leftone \intQt \dt G(\chieps) \, \thetaH
  = \iO G(\chieps(t)) \thetaH(t)
  - \iO G(\chiz) \thetaH(0)
  - \intQt G(\chieps) \, \dt\thetaH
  \non
  \\
  \leftone \geq \thetamin \iO G(\chieps(t))
  - c
  - \frac 12 \intQt |G(\chieps)|^2
  - \frac 12 \intQt |\dt\thetaH|^2
  \geq \thetamin \iO G(\chieps(t))
  - c \intQt |\chieps|^2
  - c .
  \non
\Esist
The second one is the following
\Beq
  \delta \intQt \dt\thetaeps \, \Lneps\thetaeps
  = \delta \iO \Ieps(\thetaeps(t))
  - \delta \iO \Ieps(\thetaz)
  \geq \delta \iO \Ieps(\thetaeps(t))
  - c
  \non
\Eeq
where we have set
\Beq
  \Ieps(r) := \int_1^r \Lneps s \, ds
  = \frac \eps 2 \, (r^2-1) + \int_1^r \lneps s \, ds
  \quad \hbox{for $r\in\erre$}.
\Eeq
Note that $\Ieps$ is convex and bounded from below uniformly with respect
to~$\eps$ since $\lneps1=0$.
Now, we consider the \rhs\ and deal with the non-trivial terms of~it.
We~have
\Beq
  - \intQt \dt\thetaH (\thetaeps - \thetaH)
  \leq \frac 14 \intQt |\dt\thetaH|^2
  + \intQt |\thetaeps - \thetaH|^2
  \leq \intQt |\thetaeps - \thetaH|^2
  + c .
  \non
\Eeq
Next, we~consider
\Bsist
  \leftone \intQt \nabla\Lneps\thetaeps \cdot \nabla\thetaH
  = \intQt \nabla(\Lneps\thetaeps - \Lneps\thetaH) \cdot \nabla\thetaH
  + \intQt \nabla\Lneps\thetaH \cdot \nabla\thetaH
  \non
  \\
  \leftone \leq \frac \delta 8 \intQt |\nabla(\Lneps\thetaeps - \Lneps\thetaH)|^2
  + \frac 2\delta \intQt |\nabla\thetaH|^2
  + c
  \leq \frac \delta 8 \intQt |\nabla(\Lneps\thetaeps - \Lneps\thetaH)|^2
  + \cd .
  \non
\Esist
Moreover, by \eqref{stimabetaH} and \eqref{hppi} it is easily seen~that
\Bsist
  && - \intQt \betaeps(\thetaH) (\thetaeps - \thetaH)
  + \intQt \pi(\thetaeps) (\thetaeps - \thetaH)
  \non
  \\
  && \leq c \intQt \bigl( 1 + |\thetaeps| + |\piz| \bigr) \, |\thetaeps - \thetaH|
  \leq c \intQt |\thetaeps - \thetaH|^2 + c .
  \non
\Esist
We deal with the next integral of \eqref{perprimastima}
as~follows (cf.~\eqref{stimaLnthetaH})
\Bsist
  && \delta \intQt \dt\thetaeps \, \Lneps\thetaH
  = \delta \intQt \dt(\thetaeps - \thetaH) \, \Lneps\thetaH
  + \delta \intQt \dt\thetaH \, \Lneps\thetaH
  \non
  \\
  && = \delta \iO (\thetaeps(t) - \thetaH(t)) \Lneps\thetaH(t)
  - \delta \iO (\thetaz - \thetaH(0)) \Lneps\thetaH(0)
  \non
  \\
  && \quad {}
  - \delta \intQt (\thetaeps - \thetaH) \dt\Lneps\thetaH
  + \delta \intQt \dt\thetaH \, \Lneps\thetaH
  \non
  \\
  && \leq \frac \delta 2 \iO |\thetaeps(t) - \thetaH(t)|^2
  + \intQt |\thetaeps - \thetaH|^2
  + c.
  \non
\Esist
Subsequently, by the Poincar\'e
inequality~\eqref{poincare} we~have
\Bsist
  && \delta \intQt G'(\chieps) \dt\chieps \,
       (\Lneps\thetaeps - \Lneps\thetaH)
  \non
  \\
  && \leq \frac \delta {8\cO^2}
    \intQt |\Lneps\thetaeps - \Lneps\thetaH|^2
  + 2\delta \cO^2 \intQt |G'(\chieps) \dt\chieps|^2
  \non
  \\
  && \leq \frac \delta 8 \intQt |\nabla(\Lneps\thetaeps - \Lneps\thetaH)|^2
  + 2\delta \cO^2 \sup |G'|^2 \intQt |\dt\chieps|^2.
  \non
\Esist
The next integral is easy.
We have indeed
\Beq
  - \delta \intQt \nabla \Lneps\thetaH
      \cdot \nabla(\Lneps\thetaeps - \Lneps\thetaH)
  \leq \frac \delta 8 \intQt |\nabla(\Lneps\thetaeps - \Lneps\thetaH)|^2
  + c
  \non
\Eeq
while the other two terms involving $\Lneps$ are estimated owing to the
Poincar\'e inequality once more this~way
\Bsist
  && - \delta \intQt \betaeps(\thetaH) (\Lneps\thetaeps - \Lneps\thetaH)
  + \delta \intQt \pi(\thetaeps) (\Lneps\thetaeps - \Lneps\thetaH)
  \non
  \\
  && \leq \frac \delta {8\cO^2} \intQt |\Lneps\thetaeps - \Lneps\thetaH|^2
  + 2\delta \cO^2 \intQt |\pi(\thetaeps) - \betaeps(\thetaH)|^2
  \non
  \\
  && \leq \frac \delta 8 \intQt |\nabla(\Lneps\thetaeps - \Lneps\thetaH)|^2
  + c \intQt |\pi(\thetaeps)|^2
  + c \intQt |\betaeps(\thetaH)|^2
  \non
  \\
  && \leq \frac \delta 8 \intQt |\nabla(\Lneps\thetaeps - \Lneps\thetaH)|^2
  + c \intQt |\thetaeps - \thetaH|^2
  + c .
  \non
\Esist
Next, we treat the second last integral of~\eqref{perprimastima}.
We~have
\Beq
  - \intQt F(\chieps) \dt\chieps
  \leq \frac 14 \intQt |\dt\chieps|^2
  + \intQt |F(\chieps)|^2
  \leq \frac 14 \intQt |\dt\chieps|^2
  + c \intQt |\chieps|^2
  + c .
  \non
\Eeq
Finally, we deal with the last term as~follows
\Beq
  \frac 12 \iO |\chieps(t)|^2
  = \frac 12 \iO |\chiz|^2
  + \intQt \chieps \, \dt\chieps
  \leq \frac 14 \intQt |\dt\chieps|^2
  + \intQt |\chieps|^2
  + c .
  \non
\Eeq
Hence, we can recall \eqref{perprimastima} and all the estimates we have
derived.
Forgetting some nonnegative terms on the \lhs\ we~have
\Bsist
  && \frac {1 - \delta} 2 \iO |\thetaeps(t) - \thetaH(t)|^2
  + \intQt \nabla\Lneps\thetaeps \cdot \nabla\thetaeps
  + \frac \delta 2 \intQt |\nabla(\Lneps\thetaeps - \Lneps\thetaH)|^2
  \non
  \\
  && \quad {}
  + \tonde{
      \frac 12
      - 2\delta \cO^2 \sup|G'|^2
    } \intQt |\dt\chieps|^2
  + \frac 12 \iO |\nabla\chieps(t)|^2
  + \frac 12 \iO |\chieps(t)|^2
  \non
  \\
  && \leq c \intQt |\thetaeps - \thetaH|^2
  + \intQt |\chieps|^2
  + c .
  \non
\Esist
Therefore, we can choose $\delta$ small enough and apply the Gronwall lemma.
By the Poincar\'e inequality, we conclude~that
\Bsist
  && \norma{\thetaeps-\thetaH}_{\L\infty H}
  + \norma{(\Lneps'(\thetaeps))^{1/2}\,\nabla\thetaeps}_{\L2H}
  + \norma{\Lneps\thetaeps-\Lneps\thetaH}_{\L2V}
  \non
  \\
  && \quad {}
  + \norma{\dt\chieps}_{\L2H}
  + \norma\chieps_{\L\infty V}
  \leq c .
  \non
\Esist
Then, thanks to \eqref{regthetaH} and~\eqref{stimaLnthetaH}
we obtain the basic a~priori estimate
\Bsist
  && \norma\thetaeps_{\L\infty H}
  + \norma{(\Lneps'(\thetaeps))^{1/2}\,\nabla\thetaeps}_{\L2H}
  \non
  \\
  && \quad {}
  + \norma{\Lneps\thetaeps}_{\L2V}
  + \norma\chieps_{\L\infty V\cap\H1H}
  \leq c .
  \label{primastima}
\Esist

\step First consequence

As all the terms of \eqref{secondaeps} but the Laplacian are bounded
in~$\L2H$ due to \eqref{primastima} and~\eqref{lipFG}, we immediately
derive that the Laplacian is bounded as well, whence
\Beq
  \norma\chieps_{\L2W} \leq c
  \label{stimachiHdue}
\Eeq
by the homogeneous Neumann boundary condition satisfied by $\chieps$
and the general theory of elliptic equations.

\step Second consequence

We set $\Qep:=\graffe{(x,t)\in Q:\ \thetaeps(x,t)\geq1}$ and
$\Qem:=\graffe{(x,t)\in Q:\ \thetaeps(x,t)\leq1}$.
Moreover, if $v$ is a real function on~$Q$,
$v^\pm:=\max\graffe{\pm v,0}$ are its positive and negative parts.
Then, inequalities \eqref{coerclneps}~yield
\Bsist
  && \intQ \Lneps'(\thetaeps) |\nabla\thetaeps|^2
  \geq \int_{\Qem} \Lneps'(\thetaeps) |\nabla\thetaeps|^2
  \geq \int_{\Qem} |\nabla\thetaeps|^2
  \non
  \\
  && \intQ \Lneps'(\thetaeps) |\nabla\thetaeps|^2
  \geq \int_{\Qep} \Lneps'(\thetaeps) |\nabla\thetaeps|^2
  \geq \frac 12 \int_{\Qep} \frac {|\nabla\thetaeps|^2} \thetaeps
  \non
\Esist
whence \eqref{primastima} implies~that
\Beq
  \int_{\Qem} |\nabla\thetaeps|^2 \leq c
  \aand
  \int_{\Qep} |\nabla\thetaeps^{1/2}|^2 \leq c .
  \label{daprimastima}
\Eeq
Now, we set $\Oemt:=\graffe{x\in\Omega:\ \thetaeps(t)\leq1}$ and
$\Oept:=\graffe{x\in\Omega:\ \thetaeps(t)\geq1}$ \aat\ and
use~\eqref{daprimastima} in order to estimate $\nabla\thetaeps$ in a
suitable norm.
Accounting for the first~\eqref{daprimastima}, we see~that
\Beq
  \ioT \norma{\nabla\thetaeps(t)}_{L^{4/3}(\Oemt)}^2 \, dt
  \leq c \ioT \norma{\nabla\thetaeps(t)}_{L^2(\Oemt)}^2 \, dt
  = c \int_{\Qem} |\nabla\thetaeps|^2
  \leq c .
  \label{mezzastimam}
\Eeq
On the other hand, we have
$\nabla\thetaeps=2\thetaeps^{1/2}\nabla\thetaeps^{1/2}$ a.e.\ in~$\Qep$.
Therefore, using the \Holder\ inequality, we obtain \aat
\Bsist
  && \norma{\nabla\thetaeps(t)}_{L^{4/3}(\Oept)}
  \leq 2 \norma{\thetaeps^{1/2}(t)}_{L^4(\Oept)}
    \norma{\nabla\thetaeps^{1/2}(t)}_{L^2(\Oept)}
  \non
  \\
  && \leq 2 \norma{\thetaeps(t)}_{\Lx2}^{1/2}
    \norma{\nabla\thetaeps^{1/2}(t)}_{L^2(\Oept)}
  \non
\Esist
and squaring, integrating over~$(0,T)$, and owing to
\eqref{primastima} and to the second~\eqref{daprimastima}, we derive~that
\Beq
  \ioT \norma{\nabla\thetaeps(t)}_{L^{4/3}(\Oept)}^2 \, dt
  \leq 4 \norma\thetaeps_{\L\infty H}
    \int_{\Qep} |\nabla\thetaeps^{1/2}|^2 \leq c .
  \label{mezzastimap}
\Eeq
Finally, we observe that \aat\ we have
\Bsist
  && \norma{\nabla\thetaeps(t)}_{\Lx{4/3}}^{4/3}
  = \norma{\nabla\thetaeps(t)}_{L^{4/3}(\Oemt)}^{4/3}
  + \norma{\nabla\thetaeps(t)}_{L^{4/3}(\Oept)}^{4/3}
  \quad \hbox{whence}
  \non
  \\
  && \norma{\nabla\thetaeps(t)}_{\Lx{4/3}}^2
  \leq c
    \tonde{
      \norma{\nabla\thetaeps(t)}_{L^{4/3}(\Oemt)}^2
      + \norma{\nabla\thetaeps(t)}_{L^{4/3}(\Oept)}^2
    }
  \non
\Esist
so that, using \eqref{mezzastimam} and~\eqref{mezzastimap}, we
conclude~that
\Beq
  \norma{\nabla\thetaeps}_{\L2{\Lx{4/3}}} \leq c .
  \label{quattroterzi}
\Eeq

\step Second a priori estimate

We observe~that
\Beq
  (\dt\thetaeps )\, \betaeps(\thetaeps)
  = \dt\Betaeps(\thetaeps) - \Betaepst(\thetaeps)
  \aand
  \nabla\betaeps(\thetaeps)
  = \betaepsx(\thetaeps) + \betaeps'(\thetaeps) \, \nabla\thetaeps \,.
  \non
\Eeq
Therefore, if we test \eqref{primaeps} by
$\betaeps(\thetaeps)-\betaeps(\thetaH)$, rearrange, and add the same
quantity to both sides in order to take advantage of~\eqref{belowBetaeps},
we~obtain
\Bsist
  && \iO
    \bigl(
      \Betaeps(t,\thetaeps(t)) + M \eps |\thetaeps(t) - 1|
    \bigr)
  + \intQt \nabla\Lneps\thetaeps
    \cdot \betaeps'(\thetaeps) \, \nabla\thetaeps
  + \intQt \bigl( \betaeps(\thetaeps) - \betaeps(\thetaH) \bigr)^2
  \non
  \\
  && = \iO \Betaeps(0,\thetaz)
  + \intQt \Betaepst(\thetaeps) 
  + \intQt \dt\thetaeps \, \betaeps(\thetaH)
  + M \eps \iO |\thetaeps(t) - 1|
  \non
  \\
  && \quad {}
  - \intQt \nabla\Lneps\thetaeps \cdot \betaepsx(\thetaeps)
  + \intQt \nabla\Lneps\thetaeps \cdot \nabla\betaeps(\thetaH)
  \non
  \\
  && \quad {}
  + \intQt \bigl(
             \dt G(\chieps) - \betaeps(\thetaH) + \pi(\thetaeps)
           \bigr) \, (\betaeps(\thetaeps) - \betaeps(\thetaH)).
  \label{persecondastima}
\Esist
All the terms on the \lhs\ are nonnegative, and we estimate each
integral on the \rhs, separately.
The first term is bounded due to~\eqref{hpthetaz} and~\eqref{bddBetaeps}.
The next one is treated by accounting for the second
of~\eqref{disugbetaepsxt} this~way
\Bsist
  &&\intQt \Betaepst(\thetaeps) 
  \leq c \intQt \bigl( \Betaeps(\thetaeps) + |\thetaeps| + 1 \bigr)
  \non
  \\
  && \leq c \intQt  \Betaeps(\thetaeps) + c
  \leq c \intQt \bigl( \Betaeps(\thetaeps) + M \eps |\thetaeps - 1| \bigr)
  + c
  \non
\Esist
since $\thetaeps$ has already been estimated in $\L\infty H$
by~\eqref{primastima}.
For that reason, the fourth integral of \eqref{persecondastima} is bounded
as well and we deal with the third one.
We integrate it by parts and~have
\Beq
  \intQt \dt\thetaeps \, \betaeps(\thetaH)
  = \iO \thetaeps(t) \betaeps(t,\thetaH(t))
  - \iO \thetaz \betaeps(0,\thetaH(0))
  - \intQt \thetaeps \, \dt\betaeps(\thetaH)
  \non
\Eeq
and one immediately sees that the last \rhs\ is bounded, due to
\eqref{primastima} and~\eqref{stimabetaH}.
For the same reason, the next integral of \eqref{persecondastima} that
involves $\thetaH$ is bounded as well.
Now, we treat the term continainig~$\betaepsx$.
By using the first of \eqref{disugbetaepsxt} and accounting for
\eqref{primastima} and~\eqref{stimabetaH} once more, we easily~have
\Beq
  - \intQt \nabla\Lneps\thetaeps \cdot \betaepsx(\thetaeps)
  \leq \delta \intQt \bigl( |\betaeps(\thetaeps)| + 1 \bigr)^2
  + \cd
  \leq \delta \intQt
    \bigl( \betaeps(\thetaeps) - \betaeps(\thetaH) \bigr)^2
  + \cd \,.
  \non
\Eeq
Finally, owing to \eqref{lipFG}, \eqref{stimabetaH}, and~\eqref{hppi},
we~obtain
\Bsist
  && \intQt \bigl(
           \dt G(\chieps) - \betaeps(\thetaH) + \pi(\thetaeps)
         \bigr) \, (\betaeps(\thetaeps) - \betaeps(\thetaH))
  \non
  \\
  && \leq \delta \intQt |\betaeps(\thetaeps) - \betaeps(\thetaH)|^2
  + \cd \intQt (|\dt\chieps|^2 + |\thetaeps|^2 + 1)
  \non
\Esist
and the last integral is bounded due to~\eqref{primastima}.
Therefore, if we collect \eqref{persecondastima} and the inequalities
we have proved, choose $\delta$ small enough, and apply the Gronwall lemma,
we derive~that
\Beq
  \norma{\betaeps(\thetaeps)}_{\L2H} \leq c .
  \label{secondastima}
\Eeq

\step Consequence

Estimates \eqref{primastima} and \eqref{secondastima} and our assumptions
\eqref{lipFG} and \eqref{hppi} on $G$ and $\pi$ ensure~that
\Beq
  \norma{\dt\thetaeps}_{\L2\Vp} \leq c
  \label{convdttheta}
\Eeq
just by comparison in~\eqref{primaeps}.

\Blem
\label{Lemconvqo}
Assume $z,\zn\in\LQ2$, $z>0$ \aeQ, and $\zn\to z$ \aeQ.
Moreover, let $\graffe{\eps_n}$ be a positive real sequence converging
to~$0$.
Then, $\graffe{\betaepsn(\zn)}$~converges to $\beta(z)$ \aeQ.
\Elem

\Bdim
It suffices to show that, for every $\delta\in(0,1)$, we~have
\Beq
  \betaepsn(\zn) \to \beta(z)
  \quad \hbox{almost uniformly in
    $Q^\delta := \graffe{(x,t)\in Q:\ \delta\leq z(x,t)\leq1/\delta}$}.
  \non
\Eeq
Thus, fix $\delta\in(0,1)$ and $\eta>0$.
We have to show that a subset $Q^\delta_\eta\subset Q^\delta$ exists such
that $|Q\setminus Q^\delta_\eta|\leq\eta$ and $\betaepsn(\zn)\to\beta(z)$
uniformly in~$Q^\delta_\eta$,
where $|\cpto|$ stands for the Lebesgue measure in~$\erre^4$.
By the Severini-Egorov theorem, we find $Q^\delta_\eta\subset Q^\delta$
such that $|Q\setminus Q^\delta_\eta|\leq\eta$ and $\zn\to z$ uniformly
in~$Q^\delta_\eta$, and we can prove that $\betaepsn(\zn)\to\beta(z)$
uniformly in~$Q^\delta_\eta$.
Fix $\bar n$ such~that
\Beq
  \eps_n \leq \frac \delta 2
  \aand
  \frac \delta 2 \leq \zn \leq \frac 2 \delta
  \quad \hbox{in $Q^\delta_\eta$ \ for every $n\geq\bar n$}.
  \non
\Eeq
On the other hand, we~have
\Beq
  \norma{\betaepsn(\zn) - \beta(z)}_{L^\infty(Q^\delta_\eta)}
  \leq \norma{\betaepsn(\zn) - \betaepsn(z)}_{L^\infty(Q^\delta_\eta)}
  + \norma{\betaepsn(z) - \beta(z)}_{L^\infty(Q^\delta_\eta)}.
  \label{perlemma}
\Eeq
Assume now $n\geq\bar n$.
Then, $\eps_n\leq\delta$, whence $\eps_n\leq z\leq1/\eps_n$.
Thus, $\beta(z)=\betatilde_{\eps_n}(z)$.
We infer that the last term of \eqref{perlemma} is $\leq M\eps_n$
by~\eqref{erroreregol}.
On the other hand, as $\eps_n\leq\delta/2$, we can use the second
of~\eqref{lipdelta}.
Therefore, we conclude~that
\Beq
  \norma{\betaepsn(\zn) - \beta(z)}_{L^\infty(Q^\delta_\eta)}
  \leq \cd \norma{\zn - z}_{L^\infty(Q^\delta_\eta)}
  + M \eps_n
  \non
\Eeq
and deduce that $\betaepsn(\zn)$ converges to $\beta(z)$ uniformly
in~$Q^\delta_\eta$.
\Edim

\step Conclusion of the proof

The estimates \eqref{primastima}, \eqref{stimachiHdue},
\eqref{quattroterzi}, \eqref{secondastima}, and \eqref{convdttheta} proved
in the previous steps and classical weak and weak star compacness results
ensure that suitable limit functions exist in order that the following
convergences hold (at~least for a~subsequence)
\Bsist
  & \thetaeps \to \theta
  & \quad \hbox{weakly star in $\L\infty H\cap\H1\Vp$}
  \label{convtheta}
  \\
  & \Lneps\thetaeps \to \ell
  & \quad \hbox{weakly in $\L2V$}
  \label{convLntheta}
  \\
  & \chieps \to \chi
  & \quad \hbox{weakly in $\L2W\cap\H1H$}
  \label{connvchi}
  \\
  & \nabla\thetaeps \to \nabla\theta
  & \quad \hbox{weakly in $\L2{\Lx{4/3}}$}
  \label{convquattroterzi}
  \\
  & \betaeps(\thetaeps) \to \xi
  & \quad \hbox{weakly in $\L2H$}.
  \label{convbeta}
\Esist
Note that \eqref{convtheta} and \eqref{convquattroterzi} imply that
$\thetaeps$ converges to $\theta$ weakly in $\L2{\Wx{1,4/3}}$.
Now, we observe that the Sobolev exponent $(4/3)^*$ of $\Wx{1,4/3}$ is
$12/5>2$.
Hence, $\Wx{1,4/3}$ is compactly embedded in~$H$.
On the other hand, even $W$ is compactly embedded in~$H$.
Therefore, by applying \cite[Thm.~5.1, p.~58]{Lions}
and possibly taking another subsequence, we derive~that
\Beq
  \thetaeps \to \theta
  \aand
  \chieps \to \chi
  \quad \hbox{strongly in $\L2H$ and \aeQ}.
  \label{convforte}
\Eeq
This allows us to identify all the limits of the nonlinear terms.
As far as the logarithm is concerned, we note that \eqref{convtheta} and
\eqref{convLntheta} imply that $\lneps\theta$ converges to $\ell$ weakly
in~$\L2H$.
Hence, we can conclude that $\theta>0$ and $\ell=\ln\theta$ \aeQ\
(see, e.g., \cite[Prop.~2.5, p.~27]{Brezis} for a similar result).
From $\theta>0$ \aeQ\ and \eqref{convforte} for~$\thetaeps$,
we see that we can apply Lemma~\ref{Lemconvqo} and infer that
$\xi=\beta(\theta)$ \aeQ.
Finally, the limits of the remaining nonlinear terms (i.e.,~those related
to $G$, $F'$, $G'$, and~$\pi$) can~be identified by using the convergences
a.e.\ given by \eqref{convforte} and accounting for our assumptions
\eqref{hpFG} and~\eqref{hppi}.
This concludes the proof of Theorem~\ref{Esistenza}.


\section{Boundedness}
\label{Regularity}
\setcounter{equation}{0}

In this section, we prove Theorem~\ref{Regolarita}
by estimating the $L^p$-norm of~$\theta$ (or~of a suitable function of~it)
with~a constant independent of~$p$
by using a Moser type technique.
As usual, if $z$ is either a function or a real number,
the symbol $z^+$ denotes its positive part.
Moreover, it is understood that $n$~is a positive integer and $\delta$ is a
positive parameter, say $\delta\in(0,1)$.
Furthermore, we set for convenience
\Beq
  \ustar := \ln\thetamax
  \label{defustar}
\Eeq
Finally, we assume $q\leq4$ (see~\eqref{hpq}) without loss of generality.

In the sequel, we perform two a~priori estimates.
In each of them, the use of the chain rule for time derivatives has to be
justified, and the trouble is a lack of regularity for~$\dt\theta$, which
is not known to belong to $\L2H$.
The first lemma we prove overcomes such a difficulty and uses just
the regularity we actually know for~$\theta$, namely
\Beq
  \theta \in \L\infty H \cap \H1\Vp,
  \quad \theta > 0 \quad \aeQ,
  \aand \ln\theta \in \L2V
  \label{regchain}
\Eeq
the regularity for the time derivative being 
a consequence of Theorem~\ref{Esistenza} (see Remark~\ref{Piureg}).
Actually, the last of~\eqref{regchain} does not play a special role in the lemma, 
in which a general continuous increasing function 
$ \phi' :(0,+\infty) \to \erre$ is considered.

\Blem
\label{Chainrule}
Assume $\theta\in\L\infty H\cap\H1\Vp$ and $\theta>0$ \aeQ.
Moreover, let $\phi:(0,+\infty)\to\erre$ be a convex function of class~$C^1$
and assume that $\phi'(\theta)\in\L2\Vz$. Then, if
$\Phi :\erre\to(-\infty,+\infty]$ denotes the extension 
\Beq
  \Phi(r) := \phi(r)
  \quad \hbox{if $\, r>0$},
  \quad \Phi(0) := \lim_{r\to 0^+ } \phi(r),
  \aand
  \Phi(r) := +\infty
  \quad \hbox{if $\, r<0$}
  \label{defPsi}
\Eeq
then the function $t\mapsto \iO \Phi(\theta(t))$ is absolutely 
continuous on $[0,T]$ and we have 
\Beq
  \iot \< \dt\theta(s) , \phi'(\theta(s)) > \, ds
  = \iO \Phi(\theta(t))- \iO \Phi(\theta (0) )
  \quad\hbox{for every $t\in[0,T]$}.
  \label{chainrule}
\Eeq
\Elem

\Bdim
We first observe that $\Phi$ is convex, proper, and \lsc\ in $\erre$.
In addition, we notice that
\Beq
  \phi'(u) \in \partial\Phi(u) \quad \aeO
  \quad \hbox{if} \quad
  u \in H, \quad u>0 \quad \aeO, \aand \phi'(u) \in H
  \label{csdPsi}
\Eeq
and the conjugate function $\Phi^*$ of $\Phi$ satisfies
$\partial\Phi^*=(\partial\Phi)^{-1}$.
For a.a.\ $t\in(0,T)$ we set for convenience
$v(t):=\phi'(\theta(t))$ and observe that both $\theta(t)$ and $v(t)$ lie in 
$H$, and $\theta(t)>0$ \aeO. Thus, it turns out that
$v(t)\in\partial\Phi(\theta(t))$ \aeO\ by~\eqref{csdPsi}, and consequently
\Beq
\theta (t)\in\partial\Phi^* (v(t)) \quad \aeO.
\label{incl}
\Eeq
Moreover, defining the functionals
$J:H\to(-\infty,+\infty]$ and $J_0:\Vz\to(-\infty,+\infty]$ as~follows
\Beq
  J(v) := \iO \Phi^*(v)
  \quad \hbox{if $\Phi^*(v)\in\Luno$},
  \quad J(v) := +\infty
  \quad \hbox{otherwise},
  \aand
  J_0 := J|_{\Vz}
  \label{defJ}
\Eeq
we note that (cf., e.g., \cite[Prop.~2.8, p~71]{Barbu}) 
$J$ is convex, proper, and l.s.c., and 
its subdifferential operator $\partial J : H \to 2^H$ is exactly induced 
by $\partial\Phi^*$ via the almost everywhere in $\Omega$ inclusion.  
Then, for a.a.\ $t\in(0,T)$ \eqref{incl} entails $ \theta (t)\in\partial J (v(t)) $,
that is, 
\Beq
(\theta (t), v(t) - w) + J(v(t)) \leq J(w) \quad \forall \, w \in H .
\non
\Eeq
Therefore, as $v(t) \in \Vz$ (and $v(t) \in D(\partial J) \subseteq D(J)$), 
the functional $J_0$ is proper and 
\Beq
\langle \theta (t), v(t) - w \rangle + J_0(v(t)) \leq J_0 (w) \quad \forall \, w \in \Vz 
\non
\Eeq
whence the inclusion $ \theta (t)\in\partial J_0 (v(t)) $ holds for the subdifferential 
$ \partial J_0 : \Vz \to 2^{\Vp}$ as well. Then, introducing the 
conjugate functionals and the subdifferentials
\Beq
  J^* : H \to (-\infty,+\infty],
  \quad
  J_0^* : \Vp \to (-\infty,+\infty]
  \aand
   \partial J^* : H \to 2^H,
  \quad 
  \partial J_0^* : \Vp \to 2^{\Vz}
  \non
\Eeq
and recalling that
$\partial J^*=(\partial J)^{-1}$ and $\partial J_0^*=(\partial J_0)^{-1}$,
we observe~that $ v (t)\in\partial J^* ( \theta (t)) $ and $ v (t)\in\partial J_0^* ( \theta (t)) $ for a.a.\ $t\in(0,T)$.
Therefore, being understood that $\theta$ denotes the $\Vz^* -$valued continuous 
representative, from the latter we conclude that
(see, e.g., \cite[Lemma~3.3, p.~73]{Brezis} for a similar result)  
$J_0^*(\theta)$ is absolutely continuous in $[0,T]$ and
\Beq
  \int_{0}^{t} \< \dt\theta(s) , \phi'(\theta(s)) > \, ds
  = \int_{0}^{t}  \< \dt\theta(s) , v(s) > \, ds
  = J_0^*(\theta(t)) - J_0^*(\theta(0))
  \quad \forall \, t\in[0,T].
  \label{qchainrule}
\Eeq
Moreover, the same representative $\theta$ is $H-$valued
and weakly continuous from $[0,T]$ to $H$.  
Now, as $ \Phi^{**} \equiv \Phi $ we remind that for $u, \, w \in H$ one has
\Bsist 
&&J^*(u) \displaystyle := \iO \Phi (u)
  \quad \hbox{if $\, \Phi (u)\in\Luno$},
  \quad J^*(u) := +\infty
  \  \hbox{ otherwise}\label{a-gianni}\\
&&\hbox{$ w\in\partial J^* (u) $ \ if and only if \ 
$w \in \partial \Phi (u) $ \ a.e. in $\Omega$.} 
\label{b-gianni}
\Esist
We also claim that 
\Beq
  J_0^*(u) = J^*(u)
  \quad \hbox{if $ u \in H$}. 
  \label{ugpier} 
\Eeq
For $u\in H$ we have indeed
\Beq
  J_0^*(u )
  = \sup_{w\in \Vz} \{ \< u , w > - J_0(w ) \}
  \leq  \sup_{w \in H} \{ ( u , w ) - J(w ) \}
  = J^*(u).
  \non
\Eeq
On the other hand, in view of \cite[Lemma~2.3]{DK} (or, also for related results, 
\cite[Lemma~2.4 and Section~2]{bcgg}), it turns out that  
for all $w\in D(J) $ there exixts a sequence $\{w_n\} \subset V_0 $ such that
$w_n \to w $ in H and $ J_0 (w_n) = J (w_n) \to J(w) $ as $n\to\infty$, whence 
\Beq
( u , w ) - J(w ) = \lim_{n\to \infty } \{ ( u , w_n ) - J(w_n) \} =
\lim_{n\to \infty }   \{ \< u , w_n > - J_0 (w_n) \} \leq J_0^*(u )
\non
\Eeq
and consequently $J^*(u) \leq  J_0^*(u )$ as the inequality  $( u , w ) - J(w ) \leq J_0^*(u )$
holds for all $w \in H$. 
Then, \eqref{ugpier} is proved. By combining \eqref{a-gianni} and \eqref{ugpier} 
we conclude that
\Beq
  \displaystyle J_0^*(\theta(t) ) = J^*(\theta(t)) = \iO \Phi (\theta(t) )
  \quad \forall \, t\in [0,T]. 
  \non 
\Eeq
This yields the assertion of the lemma. 
\Edim

\Brem
\label{phigrande=piccolo}
We can replace $\Phi$ by $\phi$ in the \rhs\ of \eqref{chainrule} 
in a number of cases. For instance, if we know that (the continuous representative
of) $\theta$ satisfies $\theta (t)>0$ a.e. in $\Omega$ for every $t\in [0,T]$,
then $\Phi (\theta(t)) = \phi( \theta (t) ) $ a.e. in $\Omega$ for every $t$ as well. 
Similarly, we have the same conclusion, independently of the strict positivity 
of $\theta$, whenever $\Phi(0) = +\infty $. 
\Erem

\Blem
\label{Phincoerc}
Set
\Beq
  \Phin(r)
  := \int_{\thetamax}^r
    \tonde{
      e^{2\min\graffe{n,(\ln s - \ustar)^+}} - 1
    } \, ds
  \quad \hbox{for $\, r\in(0,+\infty)$}.
  \label{defPhiepsn}
\Eeq
Then, positive constants $\alphastar$ and $\Cstar$ exist such that
\Beq
  \Phin(r)
  \geq \alphastar \, e^{3\min\graffe{n,(\ln r-\ustar)^+}}
   - \Cstar
  \label{phiepsncoerc}
\Eeq
for every $r\in(0,+\infty)$ and any positive integer~$n$.
\Elem

\Bdim
Assume first $\thetamax\leq r\leq\thetamax e^n$.
Then, we~have
\Beq
  \Phin(r)
  = \int_{\thetamax}^r \bigl( (s/\thetamax)^2 - 1 \bigr) \, ds
  = \frac \thetamax 3 \, \tonde{\frac r\thetamax}^3
  - r + \frac{2\thetamax}3
  \geq \frac \thetamax 6 \, \tonde{\frac r\thetamax}^3 -\Cstar
  \non
\Eeq
for some $\Cstar>0$, whence $\alphastar:=\thetamax/6$ works in
\eqref{phiepsncoerc} for this case.
Moreover, we can assume $\Cstar\geq\alphastar$, so that
\eqref{phiepsncoerc} holds even for $r\in(0,\thetamax)$,
since $\Phin(r)=0$ for such values of~$r$.
Finally, if $r\geq\thetamax e^n$,
we~have $r\geq r':=\thetamax e^n$
and we already now that \eqref{phiepsncoerc} holds with $r=r'$.
We deduce~that
\Beq
  \Phin(r)
  \geq \Phin(r')
  \geq \alphastar \, e^{3\min\graffe{n,(\ln r' - \ustar)^+}} - \Cstar
  = \alphastar \, e^{3n} - \Cstar
  = \alphastar \, e^{3\min\graffe{n,(\ln r - \ustar)^+}} - \Cstar.
  \non
\Eeq
This concludes the proof.
\Edim

\Blem
\label{Psincoerc}
Assume $p\in[1,+\infty)$ and~set
\Beq
  \Psin(r)
  := \int_{\thetamax}^r
    \tonde{
      \min \graffe{ n , (\ln s - \ustar)^+ }^{2p-1}
    } \, ds
  \quad \hbox{for $r\in(0,+\infty)$}.
  \label{defPsin}
\Eeq
Then, we~have
\Beq
  \Psin(r)
  \geq \frac 1 {2p}
    \, \min \graffe{ n , (\ln r - \ustar)^+ }^{2p}
    \quad \hbox{for every $r\in(0,+\infty)$}.
  \label{psincoerc}
\Eeq
\Elem

\Bdim
If $\thetamax\leq r\leq\thetamax e^n$, we~have
\Beq
  \Psin(r)
  = \int_{\ustar}^{\ln r} e^y (y-\ustar)^{2p-1} \, dy
  \geq \int_{\ustar}^{\ln r} (y-\ustar)^{2p-1} \, dy
  = \frac 1{2p} (\ln r - \ustar)^{2p}
  \non
\Eeq
and \eqref{psincoerc} follow.
If instead, $r\geq\thetamax e^n$,
we observe that $\Psin(r)\geq\Psin(r')$,
where $r':=\thetamax e^n$.
On the other hand, we already now that \eqref{psincoerc} holds with $r=r'$.
Hence, we easily conclude that the desired inequality is true
even in this case.
Finally, if $r<\thetamax$, we have $\Psin(r)=0$
and \eqref{psincoerc} trivially holds.
\Edim

Now, we start estimating.

\step First a priori estimate

We~set
\Beq
  u := \ln \theta, \quad
  \wn := \min\graffe{n,(u-\ustar)^+} ,
  \aand
  \vn := e^{2\wn} - 1.
  \label{defueps}
\Eeq
We want to use $\vn$ as a test function in~\eqref{prima} and apply
Lemma~\ref{Chainrule} with $\phi=\Phin$ given by~\eqref{defPhiepsn}.
To this aim, we note that $u\in\L2 V$ and that $\vn=\phi_n(u)$, where
$\phi_n$ is a \Lip\ continuous function.
Hence, $\vn\in\L2V$.
Moreover, $\vn$~vanishes on the boundary since $\ln\thetaG\leq\ustar$
by~\eqref{hpthetaG}.
Therefore, $\vn\in\L2\Vz$.
Furthermore, $\Phin$~is a $C^1$ convex function on~$(0,+\infty)$ and
$\vn=\Phin'(\theta)$.
Hence, we are allowed both to test \eqref{prima} by $\vn$ and to apply
Lemma~\ref{Chainrule}. Note that
$\Phin(\thetaz)=0$ since $\thetaz\leq\thetamax$ by~\eqref{hpthetaz}.
Then, by extending $\phi_n (r) $ with $0$ value for $r=0$ (cf.~\eqref{defPsi})
and setting
\Beq
  f := -\beta(\thetamax) + \pi(\theta) + \dt G(\chi)
  \label{deffeps}
\Eeq
from \eqref{chainrule} we~obtain
\Beq
  \iO \Phin(\theta(t))
  + \intQt \nabla u \cdot \nabla\vn
  + \intQt \bigl( \beta(\theta) - \beta(\thetamax) \bigr) \, \vn
  = \intQt f \, \vn .
  \label{perprimareg}
\Eeq
We treat each integral, separately.
For the first one, we~have
\Beq
  \iO \Phin(\theta(t))
  \geq \alphastar \iO e^{3\wn(t)} - \Cstar
  = \alphastar \norma{e^{\wn(t)}}_{\Lx3}^3 - c
  \non
\Eeq
thanks to Lemma~\ref{Phincoerc}.
The next term of~\eqref{perprimareg} is treated this~way
\Beq
  \intQt \nabla u \cdot \nabla\vn
  = \intQt \nabla\wn \cdot \nabla\vn
  = 2 \intQt e^{2\wn} |\nabla\wn|^2
  = 2 \intQt |\nabla e^{\wn}|^2
  \non
\Eeq
and the last term on the \lhs\ is nonnegative too.
Indeed, $\vn$~is nonnegative and
$\beta(\theta)\geq\beta(\thetamax)$ where $\vn>0$
since $\beta$ is monotone.
Thus, let us consider the \rhs.
We first notice that $f\in\LQ2$.
Therefore, using the \Holder, Sobolev, and Poincar\'e inequalities
(see~\eqref{sobolev} and~\eqref{poincare}), we~have
\Bsist
  && \intQt f \, \vn
  \leq \intQt |f| \Bigl( \bigl| e^{\wn(s)} \bigr|^2 + 1 \Bigr)
  \\
  \non
  && \leq \norma f_{\LQ1}
  + \iot \norma{f(s)}_{\Lx2}
    \norma{e^{\wn(s)}}_{\Lx3}
    \norma{e^{\wn(s)}}_{\Lx6}
  \non
  \\
  && \leq c
  + \delta \iot \norma{e^{\wn(s)}}_{\Lx6}^2 \, ds
  + \cd \iot \normaH{f(s)}^2 \norma{e^{\wn(s)}}_{\Lx3}^2 \, ds
  \non
  \\
  && \leq 2\delta \iot \norma{e^{\wn(s)} - 1}_{\Lx6}^2 \, ds
  + \cd \iot \normaH{f(s)}^2 \norma{e^{\wn(s)}}_{\Lx3}^2 \, ds
  + \cd
  \non
  \\
  && \leq 2\delta \cO^4 \intQt |\nabla e^{\wn}|^2
  + \cd \iot \normaH{f(s)}^2 \norma{e^{\wn(s)}}_{\Lx3}^2 \, ds
  + \cd
  \non
  \\
  && \leq 2\delta \cO^4 \intQt |\nabla e^{\wn}|^2
  + \cd \iot \normaH{f(s)}^2 \norma{e^{\wn(s)}}_{\Lx3}^3 \, ds
  + \cd \,.
  \non
\Esist
At this point, we collect all the above estimates, choose $\delta$ small
enough, and apply the Gronwall lemma (see, e.g., \cite[Lemma~A.4,
p.~156]{Brezis}) noting that $\normaH{f(\cpto)}^2\in L^1(0,T)$ since
$f\in\LQ2$.
We~obtain
\Beq
  \norma{e^{\wn}}_{\L\infty{\Lx3}}^3 + \intQ |\nabla e^{\wn}|^2
  \leq c.
  \label{primareg}
\Eeq

\step Consequence

From \eqref{primareg} and~\eqref{immtre}, we deduce~that $\exp(\wn)$ is
bounded in~$\LQ4$.
Hence, we can let $n$ tend to infinity and infer~that
$\exp((u-\ustar )^+)\in\LQ4$.
Now, we observe that
\Beq
  e^{(u-\ustar )^+}
  = e^{-\ustar} e^{u}
  = \theta/\thetamax
  \quad \hbox{where} \quad
  \theta > \thetamax .
  \non
\Eeq
As $\theta$ is positive, we conclude~that
\Beq
  \theta \in \LQ4.
  \label{secondareg}
\Eeq
Now, we rewrite equation~\eqref{seconda} in the~form
\Beq
  \dt\chi - \Delta\chi
  = - F'(\chi) - G'(\chi) \, \theta
  \non
\Eeq
and observe that the \rhs\ belongs to~$\LQ4$
due to \eqref{secondareg} and the \Lip\ continuity of $F$ and~$G$.
By the general theory for parabolic equations, we infer that
$\dt\chi\in\LQ4$, whence~also
\Beq
  \norma{\dt G(\chi )}_{\LQ q}
  \leq c \norma{\dt G(\chi )}_{\LQ4} < +\infty
  \label{stimadtelletre}
\Eeq
since $q\leq4$, as we have assumed at the beginning of this section.
On the other hand, we can estimate the \rhs\ of~\eqref{prima} in a better
way owing to~\eqref{secondastima}.
Indeed, recalling~\eqref{hpq}, we conclude~that
\Beq
  \norma{\pi(\theta)}_{\LQ q}
  \leq c \norma\theta_{\LQ 4} + \norma\piz_{\LQ q} < +\infty.
  \label{stimaperturb}
\Eeq

\step The Moser type procedure

Our aim is to prove an iterative estimate~for
\Beq
  w := (\ln\theta - \ustar)^+
  \label{defweps}
\Eeq
depending on the parameter $p\in[1,+\infty)$.
It is understood that the values of the constant $c$ do not depend on~$p$.
We define
\Beq
  u := \ln\theta, \quad
  \wn := \min\graffe{n,(u - \ustar)^+},
  \aand
  \vn:=\wn^{2p-1}.
  \label{ridefueps}
\Eeq
By arguing as done for the first estimate, we see that $\vn\in\L2\Vz$,
thus an admissible test function for~\eqref{prima},
and that Lemma~\ref{Chainrule} can be applied with $\Phi=\Psin$ given
by~\eqref{defPsin}, by noting that $\vn=\Psin'(\theta)$ and letting $\Psin(0)=0$.
Therefore, as $\Psin(\thetaz)=0$ by~\eqref{hpthetaz}, we~obtain
\Beq
  \iO \Psin((\theta(t))
  + \intQt \nabla u \cdot \nabla\vn
  + \intQt (\beta(\theta) - \beta(\thetamax)) \vn
  = \intQt f \vn
  \label{perterzareg}
\Eeq
where $f$ is still given by~\eqref{deffeps}.
Thanks to Lemma~\ref{Psincoerc}, we immediately derive~that
\Beq
  \iO \Psin((\theta(t))
  \geq \frac 1 {2p} \iO (\wn(t))^{2p}.
  \non
\Eeq
The next term on the \lhs\ of \eqref{perterzareg} is easily treated as
follows
\Beq
  \intQt \nabla u \cdot \nabla\vn
  = (2p-1) \intQt \wn^{2p-2} |\nabla\wn|^2
  = \frac {2p-1}{p^2} \intQt |\nabla \wn^p|^2
  \geq \frac 1p \intQt |\nabla \wn^p|^2
  \non
\Eeq
and the last one is nonnegative, since $\vn\geq0$ and
$\beta(\theta)\geq\beta(\thetamax)$ where $\vn>0$.
In order to deal with the \rhs, let us observe that $f$ belongs to~$\LQ q$
thanks to \accorpa{stimadtelletre}{stimaperturb} and \eqref{dahpbeta}.
Therefore, if $q'$~denotes the conjugate exponent of~$q$, we~have
\Beq
  \intQt f \vn
  \leq \norma f_{\LQ q} \norma\vn_{\LQ{q'}}
  \leq c \norma{\wn^{2p-1}}_{\LQ{q'}}
  = c \norma{\wn^p}_{\LQ{q'(2p-1)/p}}^{(2p-1)/p}.
  \non
\Eeq
Collecting the above estimates,
we~obtain
\Beq
  \normaH{(\wn(t))^p}^2
  + \intQt |\nabla \wn^p|^2
  \leq cp \, \norma{\wn^p}_{\LQ{q'(2p-1)/p}}^{(2p-1)/p}
  \quad \hbox{for every $t\in[0,T]$}.
  \non
\Eeq
As both terms on the \lhs\ are nonnegative, each of them satisfies the same
bound.
Therefore, owing to~\eqref{immdue}, we derive~that
\Beq
  \norma{\wn^p}_{\LQ{10/3}}^2
  \leq c \tonde{ \norma{\wn^p}_{\L\infty H}^2
    + \norma{\nabla \wn^p}_{\L2H}^2 }
  \leq cp \, \norma{\wn^p}_{\LQ{q'(2p-1)/p}}^{(2p-1)/p}.
  \non
\Eeq
At this point, we note that \eqref{defweps} and $\theta \in L^2(Q)$ trivially imply that $w\in L^r(Q) $ for every $r\in [1,+\infty)$. So we 
let $n$ tend to infinity and conclude~that
\Beq
  \norma{w^p}_{\LQ{10/3}}^2
  \leq cp \, \norma{w^p}_{\LQ{q'(2p-1)/p}}^{(2p-1)/p}
  \non
\Eeq
where $w$ is given by~\eqref{defweps}.
In other words, we~have
\Beq
  \norma w_{\LQ{10p/3}}
  \leq (c p)^{1/(2p)} \, \norma w_{\LQ{q'(2p-1)}}^{(2p-1)/(2p)}.
  \non
\Eeq
Finally, using the \Holder\ inequality and terming $|Q|$
the Lebesgue measure of~$Q$, we infer~that
\Bsist
  && \norma w_{\LQ{10p/3}}
  \leq (c p)^{1/(2p)} |Q|^{1/(4p^2q')} \,
    \norma w_{\LQ{2pq'}}^{(2p-1)/(2p)}
  \non
  \\
  && \leq \bigl( c |Q|^{1/(2pq')} \bigr)^{1/(2p)} p^{1/(2p)} \,
    \norma w_{\LQ{2pq'}}^{(2p-1)/(2p)}
  \leq (cp)^{1/(2p)} \, \norma w_{\LQ{2pq'}}^{(2p-1)/(2p)}
  \non
\Esist
since $|Q|^{1/(2pq')}$ is bounded with respect to $p\geq1$.
As we can assume the last constant $c$ to be~$\geq1$, we conclude~that
\Beq
  \norma w_{\LQ{10p/3}}
  \leq (cp)^{1/(2p)} \, \norma w_{\LQ{2pq'}}^{(2p-1)/(2p)}
  \label{terzastima}
\Eeq
with $c\geq1$.

\step Conclusion of the proof

We rewrite \eqref{terzastima} in the form
\Beq
  \norma w_{\LQ{\sigma\cdot 2pq'}}
  \leq (c p)^{1/(2p)} \, \norma w_{\LQ{2pq'}}^{(2p-1)/(2p)}
  \quad \hbox{where} \quad
  \sigma := \frac 5 {3q'}
  \label{permoser}
\Eeq
and observe that $\sigma>1$ since $q'<5/3$ by~\eqref{hpq}.
Now, we apply \eqref{permoser} to the divergent sequence $\graffe{p_k}$
defined by $p_k:=\sigma^k$ and~obtain
\Beq
  \norma w_{\LQ{2p_{k+1}q'}}
  \leq (c p_k)^{1/(2p_k)} \, \norma w_{\LQ{2p_kq'}}^{(2p_k-1)/(2p_k)}.
  \label{nuovapier}
\Eeq
Setting for convenience
\Beq
  \lkeps := \ln^+ \norma w_{\LQ{2p_kq'}}
  \non
\Eeq
and taking the positive part of the logarithm
of both sides of~\eqref{nuovapier},
we derive~that
\Beq
  \lkpueps
  \leq \frac 1{2p_k} \, \ln (c p_k)
  + \frac{2p_k-1}{2p_k} \, \lkeps
  \leq \frac 1{2p_k} \, \ln (c p_k)
  + \lkeps
  \non
\Eeq
the logarithms being nonnegative since $c\geq1$.
As this holds for every $k\geq0$, we have~that
\Beq
  \ln^+ \norma w_{\LQ{2p_kq'}}
  = \lkeps
  \leq C := \lzeps + \somma i0\infty \frac 1{2p_i} \, \ln (c p_i),
  \quad \hbox{whence} \quad
  \norma w_{\LQ{2p_kq'}} \leq e^C
  \label{finemoser}
\Eeq
by noting that the series actually converges since $p_i=\sigma^i$ with
$\sigma>1$.

At this point, we can easily conclude the proof.
Indeed, from~\eqref{finemoser}, we immediately deduce that $w\in\LQ\infty$.
Hence, coming back to~\eqref{defweps}, we derive that $\theta$ is bounded
from above.


\section{Uniqueness}
\label{Uniqueness}
\setcounter{equation}{0}

In this section, we prove Theorem~\ref{Unicita}.
The tool we use is the operator $\Rie:\Vp\to\Vz$ given by the Riesz
representation theorem, namely
\Beq
  \hbox{for $v^*\in\Vp$ and $v\in\Vz$}, \quad
  v = \Rie v^*
  \quad \hbox{means} \quad
  v^* = - \Delta v.
  \label{defRiesz}
\Eeq
We note that
\Bsist
  && \< -\Delta u , \Rie v > = \iO uv
  \quad \hbox{for every $u\in\Vz$ and $v\in H$}
  \label{rieszuv}
  \\
  && \< u^* , \Rie v^* >
  = (u^*,v^*)_*
  \quad \hbox{for every $u^*,v^*\in\Vp$}
  \label{rieszuvstar}
  \\
  && \iot \< \dt u (s) , \Rie u(s) > \, ds
  = \frac 12 \, \normaVp{u(t)}^2
  - \frac 12 \, \normaVp{u(0)}^2
  \non
  \\
  && \quad \hbox{for every $u\in\H1\Vp$ and \aat}.
  \label{rieszprimo}
\Esist
In \eqref{rieszprimo}, $\normaVp\cpto$ is the norm in $\Vp$ dual to the
norm $v\mapsto\normaH{\nabla v}$ in~$\Vz$, and the symbol on the \rhs\ of
\eqref{rieszuvstar} is the corresponding inner product.
By~the Poincar\'e inequality, such norm and product in $\Vp$ are equivalent
to the standard ones and we mainly use them for convenience.
Moreover, we recall that $R=\pi-\beta$ satisfies
assumption~\eqref{pimenobetalip}.
Finally, despite of the general rule explained at the end of
Section~\ref{MainResults}, we decide to compute all the constants we
use in our estimates with care.
In particular, we denote by $\cO$ a constant that makes the following
relations~true
\Beq
  \norma v_{\Hmuno} \leq \cO \normaVp v
  \quad \hbox{for every $v\in\Vp$}
  \aand
  \normaVp v \leq \cO \normaH v
  \quad \hbox{for every $v\in H$}
  \label{disugVp}
\Eeq
as well as the analogous of \eqref{sobolev} for~$\Lx4$.
In \eqref{disugVp}, the first norm is the standard one in~$\Vp=\Hmuno$.

To prove our uniqueness result, we have to show that any pair of
solutions $\soluz1$ and $\soluz2$ to problem \pbl\ satisfying the
regularity requirements \regsoluz\ and having the $\theta$ and $\chi$
component bounded coincide, i.e., 
$\thetau=\thetad$ and $\chiu=\chid$.
So, pick such solutions.
In particular, we can define $M$ to be the maximum of the $L^\infty$-norms
of the four functions $\thetau$, $\thetad$, $\chi_1$, and~$\chi_2$.
We note that $G$, $F'$, and $G'$ are \Lip\ continuous on~$[-M,M]$, since
they are smooth by~\eqref{hpFG}.
In the sequel, $L$~is the maximum of their \Lip\ constants on such an
interval.
Moreover, we can apply Remark~\ref{Piureg} to ($\chiu$~and) $\chid$~since
$\thetad$ is bounded, and deduce that \eqref{chipiureg} holds for~$\chid$.

Now, we write equations \eqref{prima} for both solutions and test their
difference by~$\Rie(\thetau-\thetad)$.
At the same time, we write equations \eqref{seconda} for both solutions and
test their difference by~$\mu\dt(\chiu-\chid)$, where $\mu\in(0,1)$ is a
parameter whose value will be chosen later~on.
Finally, we sum the equalities we get to each other, rearrange, and add the
same integral to both sides, for convenience.
If we use the notation $\theta:=\thetau-\thetad$ and $\chi:=\chiu-\chid$,
owing to \eqref{rieszuv}--\eqref{rieszprimo} we~obtain
\Bsist
  \leftone \frac 12 \, \normaVp{\theta(t)}^2
  + \intQt (\ln\thetau - \ln\thetad) \, \theta
  + \mu \intQt |\dt\chi|^2
  + \frac \mu 2 \iO |\nabla\chi(t)|^2
  + \frac \mu 2 \iO |\chi(t)|^2
  \non
  \\
  \leftone = \iot \bigl(
              \dt G(\chiu(s)) - \dt G(\chid(s)) , \theta(s)
            \bigr)_* \, ds
  + \int_0^t \bigl( (R(\thetau) - R(\thetad))(s) , \theta(s) \bigr)_* \, ds
  \non
  \\
  \leftone \quad {}
  + \mu \intQt \bigl( F'(\chid) - F'(\chiu) \bigr) \, \dt\chi
  + \mu \intQt \bigl( G'(\chid) \thetad - G'(\chiu) \thetau \bigr)
      \, \dt\chi
  + \frac \mu 2 \iO |\chi(t)|^2 .
  \qquad
  \label{perunicita}
\Esist
The only term on the \lhs\ that needs some treatment is the second one.
We~have
\Beq
  \intQt (\ln\thetau - \ln\thetad) \, \theta
  \geq \frac 1 M \intQt |\theta|^2
  \non
\Eeq
since $0<\thetasub_i\leq M$ for $i=1,2$.
We deal with the first term on the \rhs\ and use \eqref{chipiureg} for
$\chid$ as mentioned above. 
We~get
\Bsist
  && \iot
       \bigl(
         \dt G(\chiu(s)) - \dt G(\chid(s)) , \theta(s)
       \bigr)_* \, ds
  \non
  \\
  && = \iot
       \bigl(
         G'(\chiu(s)) \dt\chi(s)
         + (G'(\chiu(s)) - G'(\chid(s))) \dt\chid(s) ,
         \theta(s)
       \bigr)_* \, ds
  \non
  \\
  && \leq \cO \iot
    \normaH{G'(\chiu(s)) \dt\chi(s)} \normaVp{\theta(s)} \, ds
  \non
  \\
&& \quad {}
  + \cO \iot \norma{G'(\chiu(s)) - G'(\chid(s))}_{\Lx4}
    \norma{\dt\chid(s)}_{\Lx4}
    \normaVp{\theta(s)} \, ds
  \non
  \\
  && \leq \cO \, L \iot \normaH{\dt\chi(s)} \normaVp{\theta(s)} \, ds
  + \cO \, L \iot \norma{\chi(s)}_{\Lx4}
    \norma{\dt\chid(s)}_{\Lx4}
    \normaVp{\theta(s)} \, ds
  \non
  \\
  && \leq \frac \mu 8 \intQt |\dt\chi|^2
  + \frac {2 \cO^2 L^2} \mu \iot \normaVp{\theta(s)}^2 \, ds
  \non
  \\
  && \quad {}
  + \iot \norma{\chi(s)}_{\Lx4}^2 \, ds
  + \frac {\cO^2 L^2} 4
    \iot \norma{\dt\chid(s)}_{\Lx4}^2 \normaVp{\theta(s)}^2 \, ds
  \non
  \\
  && \leq \frac \mu 8 \intQt |\dt\chi|^2
  + \frac {2 \cO^2 L^2} \mu \iot \normaVp{\theta(s)}^2 \, ds
  \non
  \\
  && \quad {}
  + \cO^2 \intQt |\nabla\chi|^2
  + \cO^2 \intQt |\chi|^2
  + \frac {\cO^2 L^2} 4
    \iot \norma{\dt\chid(s)}_{\Lx4}^2 \normaVp{\theta(s)}^2 \, ds .
  \non
\Esist
Moreover, we~have
\Bsist
  && \int_0^t \bigl( (R(\thetau) - R(\thetad))(s) , \theta(s) \bigr)_* \, ds
  \leq \cO \lambdap \iot \normaH{\theta(s)} \normaVp{\theta(s)} \, ds
  \non
  \\
  && \quad {}
  \leq \frac 1 {2M} \intQt |\theta|^2
  + M \cO^2 \lambdap^2 \iot \normaVp{\theta(s)}^2 \, ds 
  \non
  \\
  && \mu \intQt \bigl( F'(\chid) - F'(\chiu) \bigr) \, \dt\chi
  \leq \frac \mu 8 \intQt |\dt\chi|^2
  + 2 \mu L^2 \intQt |\chi|^2 
  \non
  \\
  && \mu \intQt \bigl( G'(\chid) \thetad - G'(\chiu) \, \thetau \bigr)
    \, \dt\chi
  = \mu \intQt (G'(\chid) - G'(\chiu)) \, \thetau \, \dt \chi
  - \mu \intQt G'(\chid) \, \theta \, \dt\chi
  \non
  \\
  && \quad {}
  \leq \mu L M \intQt |\chi| \, |\dt\chi|
  + \mu L \intQt |\theta| \, |\dt\chi|
  \non
  \\
  && \quad {}
  \leq \frac \mu 4 \intQt |\dt\chi|^2
  + 2 \mu L^2 M^2 \intQt |\chi|^2
  + 2 \mu L^2 \intQt |\theta|^2.
  \non
\Esist
Finally, we treat the last term of \eqref{perunicita} this~way
\Beq
  \frac \mu 2 \iO |\chi(t)|^2
  = \mu \intQt \chi \, \dt\chi
  \leq \frac \mu 8 \intQt |\dt\chi|^2
  + \mu \intQt |\chi|^2 .
  \non
\Eeq
At this point, we observe that $s\mapsto\norma{\dt\chid(s)}_{\Lx4}^2$
belongs to~$L^1(0,T)$ by~\eqref{chipiureg}.
Hence, if we choose $\mu$ in order that $2\mu L^2=1/(4M)$, we see that
\eqref{perunicita} and all the estimates we have performed allow us to
apply the Gronwall lemma.
This yields $\theta=0$ and $\chi=0$ \aeQ, whence
the solutions coincide and the proof is complete.



\vspace{3truemm}

\Begin{thebibliography}{10}

\bibitem{Barbu}
V. Barbu,
``Nonlinear Semigroups and Differential Equations in Banach Spaces'',
Noordhoff International Publishing, Leyden, 1976.

\bibitem{bcgg}
 V. Barbu, P. Colli, G. Gilardi, M. Grasselli,
{\em Existence, uniqueness, and longtime behavior
for a nonlinear Volterra integrodifferential equation},
Differential Integral Equations 
{\bf 13} (2000) 1233--1262.

\bibitem{BFG} 
V. Berti, M. Fabrizio, C.Giorgi, 
{\em Well-posedness for solid-liquid phase transitions 
with a fourth-order non linearity},
Phys. D {\bf 236} (2007) 13--21.

\bibitem{BII}
E. Bonetti, 
{\em Modelling phase transitions via an entropy equation: 
long-time behaviour of the solutions}, in
``Dissipative Phase Transitions'', pp.~21--42, 
Ser. Adv. Math. Appl. Sci. {\bf 71},
World Sci. Publ., Hackensack, NJ, 2006.

\bibitem{bcf} 
E. Bonetti, P. Colli, M. Fr\'emond,
{\em A phase field model with thermal memory 
governed by the entropy balance},
Math. Models Methods Appl. Sci.
{\bf 13} (2003) 1565--1588.

\bibitem{BCFG1}
E. Bonetti, P. Colli, M. Fabrizio, G. Gilardi,
{\em Global solution to a singular integrodifferential 
system related to the entropy balance},
Nonlinear Anal. {\bf 66} (2007) 1949--1979.

\bibitem{BCFG2}
E. Bonetti, P. Colli, M. Fabrizio, G. Gilardi,
{\em Modelling and long-time behaviour for phase transitions
with entropy balance and thermal memory conductivity},
Discrete Contin. Dyn. Syst. Ser.~B
{\bf 6} (2006) 1001--1026 (electronic).

\bibitem{BFI}
E. Bonetti, M. Fr\'emond, 
{\em A phase transition model with the entropy balance},
Math. Meth. Appl. Sci.
{\bf 26} (2003) 539--556.

\bibitem{BFR}
E. Bonetti, M. Fr\'emond, E. Rocca,
{\em A new dual approach for a class of phase transitions with memory:
existence and long-time behaviour of solutions}, 
J. Math. Pures Appl. (9)
{\bf 88} (2007) 455--481.

\bibitem{BRI}
E. Bonetti, E. Rocca,
{\em Global existence and long-time time behaviour for a
singular integro-differential phase-field system},
Commun. Pure Appl. Anal.
{\bf 6} (2007) 367--387.

\bibitem{Brezis}
H. Brezis,
``Op\'erateurs Maximaux Monotones et Semi-Groupes de Contractions
dans les Espaces de Hilbert'',
North-Holland Math. Stud. {\bf 5},
North\_Holland,
Amsterdam,
1973.

\bibitem{BS} 
M. Brokate and J. Sprekels, 
``Hysteresis and Phase Transitions'',  
Appl. Math. Sci. {\bf 121}, 
Springer, 
New York,
1996.

\bibitem{cag}
G. Caginalp,
{\em An analysis of a phase field model of a free boundary},
Arch. Rational Mech. Anal.
{\bf 92} (1986) 205--245.

\bibitem{cl}
P. Colli, P. Lauren{\c{c}}ot,
{\em Weak solutions to the {P}enrose-{F}ife phase field model for a
class of admissible heat flux laws},
Phys. D
{\bf 111},(1998) 311--334.

\bibitem{cls}
P. Colli, P. Lauren{\c{c}}ot, J. Sprekels,
{\em Global solution to the {P}enrose-{F}ife phase field model with
special heat flux laws}, in
``Variations of Domain and Free-Boundary 
Problems in Solid Mechanics'', pp.~181--188, 
Solid Mech. Appl. {\bf 66},
Kluwer Acad. Publ., Dordrecht, 1999.

\bibitem{CLSS}
P. Colli, F. Luterotti, G. Schimperna, U. Stefanelli,
{\em Global existence for a class of generalized systems 
for irreversible phase changes},
NoDEA Nonlinear Differential Equations Appl.
{\bf 9} (2002) 255--276.

\bibitem{DK}
A. Damlamian, N. Kenmochi,
{\em Evolution equations generated by 
subdifferentials in the dual space
of $H^1(\Omega)$,}
Discrete Contin. Dynam. Systems,
{\bf 5} (1999) 269--278.

\bibitem{DiBen}
E. DiBenedetto,
``Degenerate Parabolic Equations'',
Springer-Verlag,
New York,
1993.

\bibitem{F1} 
M. Fabrizio, 
{\em Ginzburg-Landau equations and first and second
order phase transitions}, 
Internat. J. Engrg. Sci. 
{\bf 44} (2006) 529--539.

\bibitem{FM} 
M. Fabrizio, A. Morro, 
``Electromagnetism of Continuous Media'',
Oxford University Press, Oxford, 2003.

\bibitem{fipe}
P.C. Fife, O. Penrose,
{\em Interfacial dynamics for thermodynamically consistent
phase-field models with nonconserved order parameter},
{Electron. J. Differential Equations} {\bf 1995}, 
No.~16, approx.~49 pp.~(electronic).

\bibitem{fremond}
M. Fr\'emond,
``Non-smooth Thermomechanics'',
Springer-Verlag, Berlin, 2002.

\bibitem{germain}
P.~Germain,
``Cours de M\`ecanique des Milieux Continus. Tome I: Th\'eorie G\'en\'erale'', 
Masson et Cie, \'Editeurs, Paris, 1973.
    
\bibitem{G-L} 
V.L. Ginzburg and L.D. Landau, 
{\em On the theory of superconductivity}, 
Zh. Eksp. Teor. Fiz. 
{\bf 20} (1950) 
1064--1082.

\bibitem{gurtin}
M.E. Gurtin,
{\em Generalized Ginzburg-Landau and 
Cahn-Hilliard equations based on a microforce balance},
Phys.~D
{\bf 92} (1996)
178--192.

\bibitem{gp}
M.E. Gurtin, A.C. Pipkin,
{\em A general theory of the heat conduction with finite wave speeds},
Arch. Rational Mech. Anal
{\bf 31} (1968)
113--126.

\bibitem{Lions}
J.-L. Lions,
``Quelques M\'ethodes de R\'esolution des Probl\`emes aux Limites non
Lin\'eaires'', Dunod; Gauthier-Villars, Paris 1969.

\bibitem{moreau}
J.J.~Moreau, 
``Fonctionnelles Convexes'', 2\`eme \'edition, 
CNR \& Universit\`a di Roma Tor Vergata, Rome 2003.

\bibitem{Pf}
O. Penrose, P.C. Fife,
{\em Thermodynamically consistent models of phase-field
type for the kinetics of phase transitions},
Phys.~D
{\bf 43} (1990)
44--62.

\End{thebibliography}

\End{document}

\bye